\documentclass[twoside,reqno]{amsart}
\usepackage{amsfonts,amsmath,amsxtra, amsthm,amscd,mathrsfs}
\def\@underbar#1#2{\settowidth{\@tempdimb}{$#1#2$}\@tempdimb=0.8\@tempdimb
                   \ooalign{$#1#2$\crcr
                         \hfil\rule[-.5mm]{\@tempdimb}{.4pt}\hfil}}

\newcommand{\wt}{\operatorname{wt}}

\newcommand{\modp}{\ensuremath{\:(\text{mod }p)}}
\newcommand{\p}[1]{\ensuremath{\overline{#1}}}
\newcommand{\losemi}{{\otimes \kern -.78em \ltimes}}
\newcommand{\rosemi}{{\otimes \kern -.78em \rtimes}}

\newcommand{\Hom}{\ensuremath{\operatorname{Hom}}}

\newcommand{\Rad}{\ensuremath{\operatorname{Rad} }}
\newcommand{\Ker}{\ensuremath{\operatorname{Ker} }}

\newcommand{\dist}[1]{\operatorname{Dist}_{#1}(G)}
\newcommand{\Dist}[1]{\operatorname{Dist}(#1)}

\newcommand{\pr}{\operatorname{pr}}
\newcommand{\0}{\bar 0}
\newcommand{\1}{\bar 1}
\newcommand{\height}{\ensuremath{\operatorname{ht}}}
\newcommand{\Z}{\ensuremath{\mathbb{Z}}}
\newcommand{\Ext}{\ensuremath{\operatorname{Ext}}}
\newcommand{\g}{\ensuremath{\mathfrak{g}}}
\newcommand{\setof}[2]{\ensuremath{\left\{ #1 \: : \: #2 \right\}}}

\newcommand{\C}{\ensuremath{\mathbb{C}}}
\newcommand{\Q}{\ensuremath{\mathbb{Q}}}

\newcommand{\Op}{\ensuremath{\mathcal{O}_{p}}}

\newcommand{\Og}{\ensuremath{\mathcal{O}_{p}^{\nu}}}
\newtheorem{definition}{Definition}
\newtheorem{theorem}[definition]{Theorem}

\newtheorem{lemma}[definition]{Lemma}

\newtheorem{remark}[definition]{Remark}
\newtheorem{example}[definition]{Example}
\newtheorem{corollary}[definition]{Corollary}
\numberwithin{definition}{section}
\numberwithin{equation}{section}

\begin{document}
\title{Crystal structures arising from representations of $GL(m|n)$}
\author{Jonathan Kujawa}
\address{Dept. of Mathematics \\
           University of Toronto \\
           Toronto, ON, Canada M5S 3G3}
\email{kujawa@noether.uoregon.edu}

\issueinfo{00}
  {}
  {}
  {2003}
 \copyrightinfo{2003}
    {American Mathematical Society}
\date{\today}

\begin{abstract}
This paper provides results on the modular representation theory of the supergroup $GL(m|n).$  Working over a field of arbitrary characteristic, we prove that the explicit combinatorics of certain crystal graphs describe the representation theory of a modular analogue of the Bernstein-Gelfand-Gelfand category $\mathcal{O}$.  In particular, we obtain a linkage principle and describe the effect of certain translation functors on irreducible supermodules.  Furthermore, our approach accounts for the fact that $GL(m|n)$ has non-conjugate Borel subgroups and we show how Serganova's odd reflections give rise to canonical crystal isomorphisms.
\end{abstract}

\maketitle

\section{Introduction}\label{S:introduction} 
In 1995 Serganova computed the characters of the
finite dimensional irreducible representations of the Lie superalgebra
$\mathfrak{gl}(m|n,\C)$ \cite{serganova2,serganova1}.
Recently, Brundan
gave a more direct way to calculate these
characters \cite{B0}. He also provides for the first time a conjectural formula
for the characters of the irreducible representations belonging to the
$\mathfrak{gl}(m|n,\C)$-analogue of the Bernstein-Gelfand-Gelfand 
category $\mathcal{O}$. 
Brundan's approach relates the Grothendieck group of 
this category $\mathcal{O}$ to a certain $\mathfrak{gl}(\infty,\C)$-module.
To be more precise, let $\mathcal{V}$ denote the natural
$\mathfrak{gl}(\infty,\C)$-module and let $\mathcal{V}^{\spcheck}$ denote 
its dual. Brundan identifies 
the complexified Grothendieck group of 
category $\mathcal{O}$
with the $\mathfrak{gl}(\infty,\C)$-module
\begin{equation}\label{E:themodule}
\underbrace{\mathcal{V}^{\spcheck} \otimes\cdots\otimes
\mathcal{V}^{\spcheck}}_{n\:\text{times}} \otimes 
\underbrace{\mathcal{V} \otimes\cdots\otimes \mathcal{V}}_{m\:\text{times}}
\end{equation}
so that the Verma supermodules in category $\mathcal{O}$ correspond
to the natural monomial basis of (\ref{E:themodule}). 
Then Brundan's conjecture is that 
the irreducible supermodules in 
category $\mathcal{O}$ correspond to Lusztig's
dual canonical basis for (a certain completion of) 
the module (\ref{E:themodule}).
The subcategory $\mathcal F$
of $\mathcal O$ consisting of finite dimensional modules fits nicely
into this picture:
the Grothendieck group of the category $\mathcal F$
is identified with the submodule
$\bigwedge^{n}\mathcal{V}^{\spcheck} \otimes \bigwedge^{m}\mathcal{V}$
of (\ref{E:themodule}).

This article is concerned instead with the crystal structures (in the sense of Kashiwara) which underlie 
Brundan's conjecture. Actually, we work throughout the article
over an arbitrary field $k$ of characteristic $p$, considering a modular
analogue $\mathcal{O}_p$ of the usual category $\mathcal{O}$.
Remarkably, all the results at the level of crystals
remain true even if $p > 0$  provided one replaces the Lie algebra
$\mathfrak{gl}(\infty,\C)$ with 
the affine Kac-Moody algebra $\widehat{\mathfrak{sl}}(p,\C)$. 

When $p=0$ there are certain translation functors $E_r, F_r \:(r \in \Z)$
defined on category $\mathcal{O}$ which play a key role in
\cite{B0}; actually they already appeared in \cite{serganova2}
but in a slightly different form. At the level of
Grothendieck groups, these functors 
correspond to the usual Chevalley generators
of $\mathfrak{gl}(\infty,\C)$ acting on the module (\ref{E:themodule}). When $p>0$ one can define analogous functors $E_{r}, F_{r} \: (r \in \Z/p\Z)$ on category $\Op.$  Let $L(\lambda)$ be an irreducible module in category $\mathcal{O}$ or $\Op$.
In both cases we prove that the modules $E_r L(\lambda)$ and $F_r L(\lambda)$ are
either zero or else are self-dual indecomposable modules with irreducible
socle and cosocle isomorphic to $L(\tilde e_r^* (\lambda))$
and $L(\tilde f_r^* (\lambda))$ respectively. This gives a representation
theoretic definition
of operators $\tilde e_r^*, \tilde f_r^*$
on the set of weights that parametrizes the irreducible supermodules.
Our main result (Theorem~\ref{T:transfunctorsonirreps}) gives an explicit combinatorial description of these
operators, allowing us to verify that they are dual
to the crystal operators $\tilde e_r, \tilde f_r$
associated to Kashiwara's crystal basis of the module (\ref{E:themodule}).

Let us remark that our main result contains as a special case 
branching rules for representations of 
the supergroup $GL(m|n)$ in characteristic $p.$  They are a natural
extension
of Kleshchev's modular 
branching rules in the case of $GL(n)$; see e.g. \cite{BK2}.
Our proof is an adaptation of Kleshchev's original methods, based
on some explicit computations 
with certain lowering operators in the universal enveloping algebra.

An important feature of the Lie superalgebra 
$\mathfrak{gl}(m|n,\C)$ is that it has various different conjugacy
classes of Borel subalgebra. These may be parametrized by a sequence
$(\p{v}_1,\dots,\p{v}_{m+n})$ of parities 
$\p{v}_i \in \Z_2$, 
$m$ of which are equal to $\0$ and $n$ of which are equal to $\1$. 
Brundan only considers the Borel subalgebra corresponding to the 
parity sequence
$(\1,\cdots,\1,\0,\cdots,\0)$;
in the present article we consider all the conjugacy classes of 
Borel subalgebras.
It turns out that in the general case one must replace
the module (\ref{E:themodule}) with the module
\begin{equation}\label{E:themodule2}
\mathcal{V}^{\p{v}_1} 
\otimes \cdots \otimes \mathcal{V}^{\p{v}_{m+n}}
\end{equation}
where $\mathcal{V}^{\0}$ denotes $\mathcal{V}$ and $\mathcal{V}^{\1}$ denotes
$\mathcal{V}^{\spcheck}$.
(It appears that Brundan's conjecture itself also extends nicely 
to other Borel subalgebras, again replacing (\ref{E:themodule}) with
(\ref{E:themodule2}).)
Odd reflections, introduced by Serganova \cite{sergthesis} and, independently, by Dobrev and Petkova \cite{dobrev}, give a simple way to translate between the
parametrization of irreducible highest weight modules arising from
different choices of Borel subalgebra. In the final section of
the paper we explain how these odd reflections can also be interpreted as
canonical crystal isomorphisms between the crystals associated to
the modules (\ref{E:themodule2}) for different choices of parity sequences.

\vspace{2mm}
\noindent{\bf Acknowledgments.}
The author is grateful to Jonathan Brundan for invaluable advice during the course of this work.

\section{Preliminaries and Main Results}\label{S:conventions} 
\subsection{Supergroups}\label{SS:supergroups} Fix a ground field $k$ of characteristic $p$ (possibly $p=0$).  All objects discussed in this article (superspaces, superalgebras, supergroups, etc.) will be defined over $k.$  Recall that a superspace is a $\Z_{2}$-graded vector space and, given a superspace $V$ and a homogeneous vector $v \in V,$ we write $\p{v} \in \Z_{2}$ for the \emph{parity} (i.e. degree) of $v.$  We also recall that a commutative superalgebra $A$ is a $\Z_{2}$-graded associative algebra satisfying $ab=(-1)^{\p{a}\p{b}}ba$ for all homogeneous $a,b \in A.$  We additionally assume that $a^{2}=0$ for all $a \in A_{\1}$ when $p=2.$

Fix nonnegative integers $m$ and $n.$  Let $V$ be a superspace of even dimension $m$ and odd dimension $n;$ that is, $V$ is a $\Z_{2}$-graded vector space with $\dim_{k} V_{\0 }=m$ and $\dim_{k} V_{\1 }=n.$  The object of study in this article is the supergroup $G=GL(V).$  We define it using the language of group schemes (see \cite{J1}) as a functor from the category of commutative superalgebras to the category of groups:  for a commutative superalgebra $A$ let $G(A)$ be the group of all invertible even (i.e. grading preserving) automorphisms of the $A$-supermodule $V \otimes A.$  

\emph{Once and for all we fix a ordered homogeneous basis $v_{1}, \dotsc , v_{m+n}$ for $V.$}  Having made this choice, we can now introduce coordinates: given a commutative superalgebra $A,$ we identify $G(A)$ with the group of all invertible $(m+n) \times (m+n)$
matrices 
\begin{equation}\label{E:matrix}
(g_{i,j})_{1 \leq i,j \leq m+n }
\end{equation}
with $g_{i,j} \in A_{\p{v}_{i}+\p{v}_{j}}$ for all $1 \leq i,j \leq m+n.$  Namely, if we identify elements of $V \otimes A$ with column vectors via
\begin{equation*}
\sum_{i =1}^{m+n} v_{i}\otimes a_{i} \longleftrightarrow 
\left(\begin{matrix} a_{1} \\
\vdots \\
a_{m+n}
\end{matrix} \right),
\end{equation*}
then we can view such a matrix as an even endomorphism in the usual way by left multiplication.  Let $f: A \to B$ be a superalgebra homomorphism.  Under our identification, $G(f):G(A) \to G(B)$ corresponds to the morphism given by applying $f$ to the matrix entries.

Let $T$ be the subgroup of $G$ consisting of diagonal matrices.  More precisely, given a commutative superalgebra $A,$ $T$ is the functor given by setting $T(A)$ to be the subgroup of $G(A)$ consisting of all diagonal matrices.  Let $X(T)$ be the free abelian group 
\[
X(T)=\bigoplus_{i =1}^{m+n} \Z \varepsilon_{i}.
\]  We identify $X(T)$ with the character group of $T$ by identifying $\varepsilon_{i}$ with the function which picks out the $i$th diagonal entry of a diagonal matrix.
   We put a symmetric bilinear form on $X(T)$ by declaring that
\begin{equation}\label{bf}
(\varepsilon_i, \varepsilon_j)
= (-1)^{\p{v}_{i}} \delta_{i,j}.
\end{equation}  Observe that we have an action by the symmetric group $S_{m+n}$ on $X(T)$ given by $x \cdot \varepsilon_{i} = \varepsilon_{x i}$ for all $1 \leq i \leq m+n$ and $x \in S_{m+n}.$

Given $1 \leq t \leq m+n,$ let $V_{t}$ denote the subspace of $V$ generated by $v_{1}, \dotsc , v_{t}.$  We fix a choice of Borel subgroup $B$ of $G$ so that, for a commutative superalgebra $A,$ $B(A)$ is the stabilizer of the full flag 
\[
0 \subseteq V_{1} \otimes A \subseteq \dotsb \subseteq V_{m+n} \otimes A = V \otimes A.
\] 
Note that $B(A)$ equals the set of all upper triangular invertible matrices
of the form \eqref{E:matrix}.

The root system of $G$ is the set $\Phi = \{\varepsilon_i - \varepsilon_j : 1 \leq i, j \leq m+n , i\neq j\}$. There are even and odd roots, the 
parity of the root 
$\varepsilon_i - \varepsilon_j$ being $\p{v}_{i}  + \p{v}_{j} $.
Our choice of
Borel subgroup defines a set,
\begin{equation}\label{proots}
\Phi^+ = \{ \varepsilon_{i} - \varepsilon_{j} : 1 \leq i < j \leq m+n\},
\end{equation}
of positive roots.
The corresponding dominance order on $X(T)$ is denoted $\leq$ and is
defined by
$\lambda \leq \mu$ if $\mu - \lambda \in \mathbb{Z}_{\geq 0} \Phi^+$.

\subsection{The Superalgebra of Distributions}\label{S:algebraofdistributions}

There is an abstract notion of the superalgebra of distributions for a supergroup; see, for example, \cite[\textsection 3]{BK0}.  In this case however we can realize the superalgebra of distributions for $G,$ $\Dist{G},$ explicitly as the reduction modulo $p$ from an analogue of Kostant's $\Z$-form for the Lie superalgebra over $\C$ corresponding to $G.$   This Lie superalgebra consists of the set of $(m+n) \times (m+n)$ matrices over $\C$ with homogeneous basis given by the matrix units $e_{i,j}$ ($1 \leq i,j \leq m+n$) and with the parity of $e_{i,j}$ equal to $\p{v}_{i}+\p{v}_{j}.$  The superbracket of this Lie superalgebra is given by 
\begin{align}\label{meq}
[e_{i,j}, e_{k,l}]= \delta_{j,k}e_{i,l} - (-1)^{(\p{v}_{i} + \p{v}_{j})(\p{v}_{k}+\p{v}_{l})}\delta_{i,l}e_{k,j}.
\end{align}
Let ${U}_{\mathbb{C}}$ denote the universal enveloping superalgebra of 
this Lie superalgebra.  
By the PBW theorem for Lie superalgebras \cite{K}, 
${U}_{\mathbb{C}}$ has basis consisting of all monomials
\[
\prod_{\substack{1 \leq i,j \leq m+n\\\p{v}_{i} + \p{v}_{j}  = \0}}
e_{i,j}^{a_{i,j}}\prod_{\substack{1 \leq i,j \leq m+n\\\p{v}_{i} + \p{v}_{j} =\1}}
e_{i,j}^{d_{i,j}}
\]
where $a_{i,j} \in \mathbb{Z}_{\geq 0}$, $d_{i,j} \in \{0,1 \}$, and the product is taken in any fixed order.  
We shall write $h_{i} =e_{i,i}$ for short.

Define the \emph{Kostant $\mathbb{Z}$-form} 
${U}_{\mathbb{Z}}$ to be the $\mathbb{Z}$-subalgebra of 
${U}_{\mathbb{C}}$ generated by elements
$e_{i,j}\:(1 \leq i,j \leq m+n, \p{v}_{i}  + \p{v}_{j}  = \1)$,
$e_{i,j}^{(r)}\: 
(1 \leq i,j \leq m+n, i \neq j, \p{v}_{i}  + \p{v}_{j}  = \0, r \geq 1),$
and
$ \binom{h_i}{r}\:(1 \leq i \leq m+n, r \geq 1)$.
Here, $e_{i,j}^{(r)} := e_{i,j}^r / (r!)$ and 
$\binom{h_{i}}{r} := {h_{i}(h_{i}-1) \dotsb (h_{i}-r+1)}/{(r!)}$.
Following the proof of \cite[Th.2]{S1}, one verifies the following:

\begin{lemma}\label{stbas}
The superalgebra ${U}_{\mathbb{Z}}$ is a free
$\mathbb{Z}$-module with basis given by the set of all monomials
of the form
$$
\prod_{\substack{1 \leq i,j \leq m+n  \\
                   i \neq j, \p{v}_{i}+\p{v}_{j} = \0}} e_{i,j}^{(a_{i,j})}
\prod_{1 \leq i \leq m+n} \binom{h_{i}}{r_{i}}
\prod_{\substack{ 1 \leq i,j \leq m+n\\
\p{v}_{i}+\p{v}_{j} = \1}} e_{i,j}^{d_{i,j}}
$$
for all $a_{i,j}, r_{i} \in \mathbb{Z}_{\geq 0}$ and $d_{i,j} 
\in \{0,1 \}$, where the product is taken in any fixed order.
\end{lemma}

The enveloping superalgebra ${U}_{\mathbb{C}}$ is a 
Hopf superalgebra in the canonical way and, furthermore, this structure restricts to make 
${U}_{\mathbb{Z}}$ a Hopf superalgebra over $\mathbb{Z}$.
Finally, set
\[
\Dist{G}= k \otimes_{\mathbb{Z}} {U}_{\mathbb{Z}},
\]
naturally a Hopf superalgebra over $k$.
We will abuse notation by using the same symbols
$e_{i,j}^{(r)}, \binom{h_i}{r}$ etc... for the
canonical images of these elements of ${U}_{\mathbb{Z}}$ in
$\Dist{G}$.

It is also easy to describe the superalgebras of distributions of
our various natural subgroups of $G$ as subalgebras
of $\dist{}$.
For example, $\Dist{T}$ is the subalgebra generated by all
$ \binom{h_i}{r}\:( 1 \leq i \leq m+n, r \geq 1)$,
and 
$\Dist{B}$ is the subalgebra generated by $\Dist{T}$ and all $e^{(a)}_{i, j}$ where $1 \leq i,j \leq m+n$ with $i <j,$ and $a \in \Z_{\geq 0} \text{ if } \p{v}_{i}+\p{v}_{j}=\0\text{ and } a \in \{0,1 \} \text{ if } \p{v}_{i}+\p{v}_{j}=\1.$

Let us describe the category of $\Dist{G}$-supermodules.  The objects are all left $\Dist{G}$-modules which are $\Z_{2}$-graded; i.e., $k$-superspaces, $M,$ satisfying $\Dist{G}_{r}M_{s} \subseteq M_{r+s}$ for all $r, s \in \Z_{2}.$   A morphism of $\Dist{G}$-supermodules is a linear map $f: M\to M'$ satisfying $f(xm)=(-1)^{\p{f}\;\p{x}}xf(m)$ for all $m \in M$ and all $x \in \Dist{G}.$  Note that this definition makes sense as stated only for homogeneous elements; it should be interpreted via linearity in the general case.  We emphasize that we allow \emph{all} morphisms and not just graded (i.e. \emph{even}) morphisms.   However, we note that for superspaces $M$ and $M'$ the space $\Hom_{k}(M,M')$ is naturally $\Z_{2}$-graded by $f \in \Hom_{k}(M,M')_{r}$ if $f(M_{s}) \subseteq M'_{s+r}$ for all $r, s \in \Z_{2}.$  This gives a $\Z_{2}$-grading on $\Hom_{\Dist{G}}(M,M') \subseteq \Hom_{k}(M,M').$  The category of $\Dist{G}$-supermodules is \emph{not} an abelian category.  However, the \emph{underlying even category}, consisting of the same objects but only the even morphisms, is an abelian category.  This, along with the parity change functor, $\Pi,$ which simply interchanges the $\Z_{2}$-grading of a supermodule, allows us to make use of the tools of homological algebra.

For $\lambda=\sum_{i =1}^{m+n} \lambda_{i}\varepsilon_{i} \in X(T)$  
and a $\dist{}$-supermodule $M$, define the {\em $\lambda$-weight space of $M$}
to be
\begin{equation}\label{wtsp}
M_{\lambda}=\left\{m \in M\: : \:
\binom{h_{i}}{r}m=\binom{\lambda_{i}}{r}m \mbox{ for all } 1 \leq i \leq m+n,
r \geq 1 \right\}. 
\end{equation}
We call a $\dist{}$-supermodule $M$ {\em integrable}
if it is locally finite over $\dist{}$ and satisfies
$M = \sum_{\lambda \in X(T)} M_{\lambda}$.   The category of $G$-supermodules can naturally be identified with the category of integrable $\Dist{G}$-supermodules \cite[Corollary 3.5]{BK0}.

\subsection{Highest Weight Theory}\label{SS:serg}

Given $\lambda \in X(T)$ let 
\[
D(\lambda)=\setof{\mu\in X(T)}{\mu \leq \lambda \text{ in the dominance order}}.
\]  Let $\Op$ denote the the full subcategory of the category of all $\Dist{G}$-supermodules consisting of supermodules $M$ 
such that $M = \bigoplus_{\lambda \in X(T)} M_\lambda,$ $\dim_{k} M_{\lambda} < \infty$ for all $\lambda \in X(T),$ and there are $\lambda^{(1)}, \dotsc , \lambda^{(r)} \in X(T)$ so that $M_{\mu} \neq 0$ implies $\mu \in \bigcup_{i=1}^{r}D(\lambda^{(i)})$ for any $\mu \in X(T).$   Note that any supermodule in $\Op$ is locally finite over $\Dist{B}.$  We also remark that in the case $p=0$ Brundan's category $\mathcal{O}$ discussed in the introduction is a full subcategory of $\mathcal{O}_{0}.$  From now on we will assume all $\Dist{G}$-supermodules under discussion are objects in $\Op.$

For $\lambda \in X(T)$, we have the {\em Verma module}
$$
M(\lambda) := \Dist{G} \otimes_{\Dist{B}} k_\lambda,
$$
where $k_\lambda$ denotes $k$ viewed as a 
$\Dist{B}$-supermodule of weight $\lambda$ concentrated in degree $\0.$  Note that by Lemma~\ref{stbas} it follows that $M(\lambda)$ is an object in $\Op.$
We say that a homogeneous vector $v$ in a $\Dist{G}$-supermodule $M$
is a \emph{primitive vector of weight $\lambda$} (or simply a \emph{primitive vector}) if 
$\Dist{B} v \cong k_\lambda$ as a $\Dist{B}$-supermodule.
Familiar arguments exactly as for semisimple Lie algebras
over $\mathbb C$ show:

\begin{lemma}\label{hw} Let $\lambda \in X(T)$.
\begin{itemize}
\item[(i)]
The $\lambda$-weight space of $M(\lambda)$ is $1$-dimensional,
and all other weights of $M(\lambda)$ are $< \lambda$ in the dominance order.
\item[(ii)]
Any non-zero quotient of $M(\lambda)$
is generated by a primitive vector of weight $\lambda$,
unique up to scalars.
\item[(iii)]
Any $\dist{}$-supermodule generated by a primitive
vector of weight $\lambda$ is isomorphic to a quotient of 
$M(\lambda)$.
\item[(iv)]
$M(\lambda)$ has a unique maximal submodule $\Rad M(\lambda)$ and, hence, an irreducible quotient $L(\lambda):= M(\lambda)/\Rad M(\lambda).$
The $\{L(\lambda)\}_{\lambda \in X(T)}$ give
a complete set of pairwise non-isomorphic irreducibles in
$\Op$.
\end{itemize}
\end{lemma}

In this way, we get a parametrization of the irreducible objects
in $\Op$ by their highest weights with respect to the ordering
$\leq$.  

Given a $\Dist{G}$-supermodule, $M,$ we can consider its graded dual 
\[
M^{*}:=\bigoplus_{\lambda \in X(T)} \Hom_{k}(M_{\lambda}, k)
\]
with the usual $\Z_{2}$-grading and $\Dist{G}$ action.  We have an automorphism, $\tau,$ of $\Dist{G}$ induced by $e_{i,j} \mapsto - (-1)^{\p{v}_{i}(\p{v}_{i}+\p{v}_{j})}e_{j,i}$ (the negative of the \emph{supertranspose}).  Twisting $M^{*}$ by $\tau$ yields a new $\Dist{G}$-supermodule, which we call the \emph{contravariant dual} and denote by $M^{\tau}.$  In particular, for $\lambda \in X(T)$ we define 
\begin{align}\label{E:defofcovermaandcostandard}
W(\lambda)&=M(\lambda)^{\tau},
\end{align} the \emph{co-Verma supermodule of highest weight $\lambda.$}  We remark that for a module $M$ in category $\Op$ we have $(M^{\tau})^{\tau} \cong M$ and the characters of $M$ and $M^{\tau}$ coincide. In particular, we have $L(\lambda)^{\tau} \cong L(\lambda)$ for any $\lambda \in X(T).$

\subsection{Crystals}\label{SS:crystals}Let us recall the general definition of a crystal in the sense of Kashiwara \cite[7.2]{Kash}.  Assume we have the following data:  
\[
\begin{array}{l}
P = \text{ a free $\Z$-module (called the weight lattice)}\\
I = \text{ an index set}\\
\alpha_{i} \in P  \text{ for all $i \in I$ (called a simple root)}\\
h_{i} \in P^{*}=\Hom_{\Z}(P, \Z)  \text{ for all $i \in I$ (called a simple coroot)}\\
\langle -,- \rangle: P \times P \to \Q \text{ a symmetric bilinear form.}
\end{array}
\]
Additionally, we assume the data satisfies the following axioms:
\[
\begin{array}{l}
\langle \alpha_{i},\alpha_{i} \rangle \in 2\Z_{> 0} \text{ for all $i \in I$}\\
h_{i}(\lambda)=\frac{2\langle \alpha_{i}, \lambda \rangle}{\langle \alpha_{i}, \alpha_{i}\rangle} \text{ for all $i \in I$}\\
\langle \alpha_{i}, \alpha_{j} \rangle \leq 0 \text{ for all $i, j \in I$ with $i \neq j.$}
\end{array}
\]

 With this fixed data, we define a \emph{crystal} $\mathcal{B}$ as a set along with maps 
\begin{align*} \tilde{e}_{i}, \tilde{f}_{i}&: \mathcal{B} \to \mathcal{B} \sqcup \{0 \} \text{ (for $i \in I$)}\\
               \varepsilon_{i}, \varphi_{i}&: \mathcal{B} \to \Z \sqcup \{ -\infty \} \text{ (for $i \in I$)}\\
               \wt &: \mathcal{B} \to P
\end{align*} subject to the following axioms:
\begin{enumerate}
\item [(C1)] $\varphi_{i}(b)=\varepsilon_{i}(b)+\frac{2\langle{\alpha_{i}, \wt(b)\rangle}}{\langle \alpha_{i}, \alpha_{i}\rangle};$
\item [(C2)] if $\tilde{e}_{i}(b)\neq 0,$ then $\varepsilon_{i}(\tilde{e}_{i}(b))=\varepsilon_{i}(b)-1,$ $\varphi_{i}(\tilde{e}_{i}(b))=\varphi_{i}(b)+1,$ and $\wt(\tilde{e}_{i}(b))=\wt(b)+\alpha_{i};$
\item [(C3)] if $\tilde{f}_{i}(b)\neq 0,$ then $\varepsilon_{i}(\tilde{f}_{i}(b))=\varepsilon_{i}(b)+1,$ $\varphi_{i}(\tilde{f}_{i}(b))=\varphi_{i}(b)-1,$ and $\wt(\tilde{f}_{i}(b))=\wt(b)-\alpha_{i};$
\item [(C4)] $b_{1}=\tilde{f}_{i}(b_{2})$ if and only if $\tilde{e}_{i}(b_{1})=b_{2};$
\item [(C5)] if $\varphi_{i}(b)=-\infty,$ then $\tilde{e}_{i}(b)=\tilde{f}_{i}(b)=0;$ 
\end{enumerate} where $b, b_{1}, b_{2} \in \mathcal{B}$ and $i\in I.$

We also remind the reader of the notion of the tensor product of two crystals.  If $\mathcal{B}_{1}$ and $\mathcal{B}_{2}$ are crystals, then set $\mathcal{B}_{1}\otimes \mathcal{B}_{2}=\{b_{1}\otimes b_{2}: b_{1}\in \mathcal{B}_{1}, b_{2}\in \mathcal{B}_{2} \}$ and
\begin{align*}
\wt(b_{1}\otimes b_{2})&=\wt(b_{1}) + \wt(b_{2}),\\
\varepsilon_{i}(b_{1}\otimes b_{2}) &= \max\left( \varepsilon_{i}(b_{1}), \varepsilon_{i}(b_{2})-\frac{2\langle\alpha_{i}, \wt(b_{1})\rangle}{\langle \alpha_{i},\alpha_{i} \rangle}\right),\\
\varphi_{i}(b_{1}\otimes b_{2}) &= \max\left(\varphi_{i}(b_{2}), \varphi_{i}(b_{1})+\frac{2\langle\alpha_{i}, \wt(b_{2})\rangle}{\langle \alpha_{i},\alpha_{i}\rangle}\right),\\
\tilde{e}_{i}(b_{1}\otimes b_{2})&=\begin{cases} \tilde{e}_{i}b_{1} \otimes b_{2}, &\text{ if $\varphi_{i}(b_{1}) \geq \varepsilon_{i}(b_{2})$};\\
                                               b_{1} \otimes \tilde{e}_{i}b_{2}, &\text{ if $\varphi_{i}(b_{1}) < \varepsilon_{i}(b_{2})$};\\
                                    \end{cases}\\
\tilde{f}_{i}(b_{1}\otimes b_{2})&=\begin{cases} \tilde{f}_{i}b_{1} \otimes b_{2}, &\text{ if $\varphi_{i}(b_{1}) > \varepsilon_{i}(b_{2})$};\\
                                               b_{1} \otimes \tilde{f}_{i}b_{2}, &\text{ if $\varphi_{i}(b_{1}) \leq \varepsilon_{i}(b_{2})$}.\\
                                    \end{cases}
\end{align*}  If $M_{1}$ and $M_{2}$ are Lie algebra modules with associated crystals $\mathcal{B}_{1}$ and $\mathcal{B}_{2},$ respectively, then by \cite[Thm. 4.1]{Kash} the crystal associated to $M_{1}\otimes M_{2}$ is $\mathcal{B}_{1}\otimes \mathcal{B}_{2}.$

\subsection{Affine Lie Algebras}\label{SS:Liealgebras}

Recall that we have a fixed ground field $k$ of characteristic $p$.  There are two cases to consider: when $p =0$ and when $p >0.$  In each case we define the requisite Cartan datum in the notation of subsection~\ref{SS:crystals} and use this data to define an affine Lie algebra over $\C$, $\g,$ in the manner of \cite{Kac2}. 

We first consider the case when $p =0.$  Let $P=\sum_{r \in \Z} \Z \gamma_{r}.$ The index set is $\Z/p\Z=\Z.$  Define the \emph{simple roots} by $\alpha_{r}=\gamma_{r}-\gamma_{r+1}$ for $r \in \Z.$   We define a nondegenerate symmetric bilinear form $\langle -, - \rangle: P \times P \to\Q$ by setting $\langle \gamma_{r},\gamma_{s} \rangle = \delta_{r,s}$ for $r, s \in \Z.$  Observe that 
\[
\langle \alpha_{r}, \alpha_{s} \rangle = \begin{cases} 2, &\text{ if $r=s;$}\\
                                                       -1, &\text{ if $r=s \pm 1;$}\\
                                                       0, &\text{ otherwise;} 
\end{cases}
\] for all $r,s \in \Z.$
 Using the form we identify $P$ and $P^{*}$ via $x \leftrightarrow \langle x,- \rangle.$  Under this identification the simple coroot $h_{r}$ is $\alpha_{r}$ for all $r \in \Z.$

Now we consider the case when $p >0.$  Then we let $P=\Z\delta \oplus \bigoplus_{r \in \Z/p\Z} \Z \Lambda_{r}$.  Again, the index set is $\Z/p\Z.$  The \emph{simple roots} are defined by $\alpha_{r}=2\Lambda_{r}-\Lambda_{r-1}-\Lambda_{r+1} + \delta_{r,0}\delta$ for $r \in \Z/p\Z.$  Let $\langle - , - \rangle: P \times P \to \Q$ be the nondegenerate bilinear form determined by requiring $\delta, \Lambda_{0}, \dots , \Lambda_{p-1}$ and $\Lambda_{0}, \alpha_{0}, \dots , \alpha_{p-1}$ to be dual bases with respect to the form. Observe that if $p>2,$ then
\[
\langle \alpha_{r}, \alpha_{s} \rangle = \begin{cases} 2, &\text{ if $r\equiv s \modp;$}\\
                                                       -1, &\text{ if $r\equiv s \pm 1 \modp;$}\\
                                                       0, &\text{ otherwise;} 
\end{cases}
\] and if $p=2,$ then
\[
\langle \alpha_{r}, \alpha_{s} \rangle = \begin{cases} 2, &\text{ if $r\equiv s \:(\text{mod } 2);$}\\
                                                       -2, &\text{ if $r\equiv s + 1 \:(\text{mod } 2);$}
\end{cases}
\] for all $r,s \in \Z/p\Z.$ In particular this implies the form is symmetric.  Using the form we identify $Q:=\Z\Lambda_{0} \oplus \bigoplus_{r \in \Z/p\Z} \Z \alpha_{r} \subseteq P$ with $P^{*}$ via $x \leftrightarrow \langle x, - \rangle.$   Under this identification, the simple coroot $h_{r}$ is $\alpha_{r}$ for all $r \in \Z/p\Z.$  Finally, given $a \in \Z$ we write $a = pd+s$ with $d \in \Z$ and $s=1, \dots , p$ and define $\gamma_{a} \in P$ by
\[
\gamma_{a}=   \Lambda_{s}-\Lambda_{s-1}-d\delta.
\] Observe that if $a=pd+s$, then $\gamma_{a}-\gamma_{a+1}=\alpha_{s}.$

In both cases we define a Lie algebra over $\C,$ $\g,$ generated by $\mathfrak{h}:=P \otimes_{\Z} \C$ and $\{E_{r},F_{r}: r \in \Z/p\Z \}$ subject to the relations
\begin{align*}
[E_{r},F_{s}]&=\delta_{r,s}\alpha_{r}\\
[H,H']&=0\\
[H,E_{r}]&=\langle \alpha_{r}, H \rangle E_{r}\\
[H,F_{r}]&=-\langle \alpha_{r},H \rangle F_{r}\\
(\operatorname{ad} E_{r})&^{1-\langle\alpha_{r},\alpha_{s}\rangle} (E_{s})=0\\
(\operatorname{ad} F_{r})&^{1-\langle\alpha_{r},\alpha_{s}\rangle} (F_{s})=0
\end{align*} for all $r,s \in \Z/p\Z,$ all $H, H' \in \mathfrak{h},$ and where $\langle -,- \rangle$ denotes the bilinear form on $P$ extended to $\mathfrak{h}.$   Note that if $p =0,$
then $\g=\mathfrak{gl}_{\infty}(\C).$  If $p>0,$ then $\g=\widehat{\mathfrak{sl}}_{p}(\C).$

\subsection{The Crystal $\mathcal{B}$}\label{SS:crystalBmn}

We are now prepared to describe the crystal which plays a central role in this paper. Let $\mathcal{V}$ denote the natural ``evaluation'' $\g$-module with basis $\{x_{b}: b \in \Z \}$ and action given by
\begin{align*}
E_{r}x_{b}&=\begin{cases} x_{b-1}, &\text{ if $r+1\equiv b \modp;$}\\
                         0, &\text{ otherwise;}
\end{cases}\\
F_{r}x_{b}&=\begin{cases} x_{b+1}, &\text{ if $r\equiv b \modp;$}\\
                         0, &\text{ otherwise;}
\end{cases}\\
Hx_{b}&=\langle H, \gamma_{b}\rangle x_{b} \text{ for all $H \in \mathfrak{h},$}
\end{align*} where $\g$ is the affine Lie algebra defined in subsection~\ref{SS:Liealgebras}.  We say a vector $ x \in \mathcal{V}$ is of weight $\lambda \in P$ if $Hx=\langle H, \lambda\rangle x$ for all $H \in \mathfrak{h}.$  

 There is a crystal $(\mathcal{B}_{\0}, \tilde{e}_{r}, \tilde{f}_{r}, \varepsilon_{r}, \varphi_{r}, \wt)$ associated to the module $\mathcal{V}$ (for both $p =0$ and  $p > 0$) where the underlying set $\mathcal{B}_{\0}$ is $\{x_{b}: b\in \Z \}$ (the given basis).  The crystal operators are defined by $\tilde{e}_{r}=E_{r}$ and $\tilde{f}_{r}=F_{r}$ for all $r \in \Z/p\Z$ and, given $b \in \Z$, $\varepsilon_{r}(x_{b})=1,$ if $r+1 \equiv b \modp$ and is zero otherwise; and $\varphi_{r}(x_{b})=1,$ if $r \equiv b \modp$ and is zero otherwise.  Finally, $\wt$ is the usual weight function on $\mathcal{V},$ hence $\wt (x_{b})=\gamma_{b}$ for all $b \in \Z.$  We leave it to the reader to verify the crystal axioms. 

We have an automorphism of $\g$ given by 
\begin{align*}
E_{r} &\mapsto F_{r},\\
F_{r} &\mapsto E_{r},\\
H &\mapsto -H,
\end{align*} for all $r\in \Z/p\Z$ and all $H \in \mathfrak{h}.$  We can twist $\mathcal{V}$ by this automorphism and obtain a new $\g$-module, $\mathcal{V}^{\spcheck{}}.$  This module also has an associated crystal $(\mathcal{B}_{\1}, \tilde{e}_{r}, \tilde{f}_{r}, \varepsilon_{r}, \varphi_{r}, \wt)$ which is, roughly speaking, $\mathcal{B}_{\0}$ with the roles of $\tilde{e}_{r}$ and $\tilde{f}_{r}$ interchanged.  Namely, the crystal $\mathcal{B}_{\1}$ is the set $\{x_{b}^{\spcheck{}}: b \in \Z \}$ (the given basis of $\mathcal{V}$ viewed as a basis of $\mathcal{V}^{\spcheck}$),  where we set $\tilde{e}_{r}(x_{b}^{\spcheck{}})=\tilde{f}_{r}(x_{b})^{\spcheck{}},$ $\tilde{f}_{r}(x_{b}^{\spcheck{}})=\tilde{e}_{r}(x_{b})^{\spcheck{}},$ $\varepsilon_{r}(x_{b}^{\spcheck{}})=\varphi_{r}(x_{b}),$ $\varphi_{r}(x_{b}^{\spcheck{}})=\varepsilon_{r}(x_{b}),$ and $\wt(x_{b}^{\spcheck{}})=-\wt(x_{b}),$ for all $r\in \Z/p\Z$ and all $b \in \Z.$  

The tensor product of crystals 
\[
\left( \mathcal{B}_{\p{v}_{1}} \otimes \dotsb \otimes \mathcal{B}_{\p{v}_{m+n}}, \tilde{e}_{r}, \tilde{f}_{r}, \varepsilon_{r}, \varphi_{r}, \wt\right)
\] is then the crystal associated to the $\g$-module $\mathcal{V}^{\p{v}_{1}}\otimes\dotsb \otimes \mathcal{V}^{\p{v}_{m+n}}$ where $\mathcal{V}^{\0}:=\mathcal{V}$ and $\mathcal{V}^{\1}:=\mathcal{V}^{\spcheck{}}.$ 
The key combinatorial object in this article is a different crystal structure $(\mathcal{B}, \tilde{e}_{r}^{*}, \tilde{f}_{r}^{*}, \varepsilon_{r}^{*}, \varphi_{r}^{*}, \wt)$ on the same underlying set,
\[
\mathcal{B}:=\setof{x_{b_{1}}^{\p{v}_{1}} \otimes \dots  \otimes x_{b_{m+n}}^{\p{v}_{m+n}}}{b_{1}, \dots, b_{m+n} \in \Z},
\] where $x_{b}^{\0}:=x_{b} \in \mathcal{B}_{\0}$ and $x_{b}^{\1}:=x_{b}^{\spcheck} \in \mathcal{B}_{\1}.$  We call this the \emph{dual crystal structure} following Brundan \cite{B0}.  The dual crystal operators are defined by 
\[\begin{array}{cc}
\tilde{e}^{*}_{r}(x):=-\tilde{f}_{-1-r}(-x), & \varepsilon_{r}^{*}(x):=\varphi_{-1-r}(-x),\\
\tilde{f}^{*}_{r}(x):=-\tilde{e}_{-1-r}(-x), & \varphi_{r}^{*}(x):=\varepsilon_{-1-r}(-x),
\end{array}
\] where for 
\[
x=x_{b_{1}}^{\p{v}_{1}} \otimes \dots \otimes x_{b_{m+n}}^{\p{v}_{m+n}} \in \mathcal{B}
\]
we define
\[
-x=x_{-b_{1}}^{\p{v}_{1}} \otimes \dots  \otimes x_{-b_{m+n}}^{\p{v}_{m+n}}.
\] The weight function is given by 
\[
\wt\left( x_{b_{1}}^{\p{v}_{1}} \otimes \dots  \otimes x_{b_{m+n}}^{\p{v}_{m+n}}\right)= \sum_{i=1}^{m+n} (-1)^{\p{v}_{i}}\gamma_{b_{i}}.
\]  See \cite{B0} for a discussion of the sense in which these crystal structures are dual to one another.  In the next subsection we will give a more explicit combinatorial description of the dual crystal $\mathcal{B}.$

\subsection{The Crystal Structure on $X(T)$}\label{SS:crystalofXofT}
We now lift the dual crystal structure on $\mathcal{B}$ to $X(T).$   To do so we require some additional notation.  Let $\rho \in X(T)$ denote the unique element which satisfies the following conditions:
$$(\rho, \varepsilon_{m+n} ) = \begin{cases} 1, & \text{ if } \p{v}_{m+n}=\0;\\
                                                               0, &  \text{ if } \p{v}_{m+n}=\1;
\end{cases}$$
and for $1 \leq i \leq m+n-1,$
\[
 ( \rho, \varepsilon_{i}-\varepsilon_{i+1}) = \begin{cases} 1, & \text{ if } \p{v}_{i}=\p{v}_{i+1}=\0 ;\\
                                                                             -1, & \text{ if } \p{v}_{i}=\p{v}_{i+1}=\1 ;\\
                                                                             0, & \text{ if } \p{v}_{i} \neq \p{v}_{i+1}.
\end{cases}
\]
For $1 \leq j \leq m+n,$ let
\begin{equation}\label{E:thetadef}
\vartheta_{j}=\sum_{i=j+1}^{m+n}(-1)^{\p{v}_{i}+\p{v}_{j}},
\end{equation}  and let
\begin{equation}\label{E:bigthetadef}
\vartheta= \sum_{j =1}^{m+n}\vartheta_{j}\varepsilon_{j} \in X(T).
\end{equation}  Then observe that we have 
\begin{align}\label{E:rhodef}
\rho&= \vartheta + \sum_{\substack{1 \leq i \leq m+n\\ \p{v}_{i}=\0}} \varepsilon_{i}.
\end{align} As an example, say our fixed homogeneous basis for $V$ satisfies $\p{v}_{i}=\1 $ if $1 \leq i \leq n$ and $\p{v}_{i}=\0$ if $n+1 \leq i \leq m+n$ (the parity choice in \cite{B0}), then we have
\begin{align*}\label{E:bigthetadef2}
\vartheta&= -(m-n+1)\varepsilon_{1} +\dots +-m\varepsilon_{m}+(m-1)\varepsilon_{n+1}+\dots +\varepsilon_{m+n-1}, \\
\rho&=-(m-n+1)\varepsilon_{1} +\dots +-m\varepsilon_{m}+m\varepsilon_{n+1}+\dots +\varepsilon_{m+n}.
\end{align*}

We define a bijection $X(T) \to \mathcal{B}$ by 
\begin{equation}\label{E:crystalbijectiondef}
\lambda \mapsto x_{(\lambda + \rho, \varepsilon_{1})}^{\p{v}_{1}} \otimes \dots \otimes  x_{(\lambda+\rho, \varepsilon_{m+n})}^{\p{v}_{m+n}}.
\end{equation}
Using this bijection we can lift the dual crystal structure on $\mathcal{B}$ to the set $X(T)$.  Let us describe the combinatorics of this crystal structure explicitly.  It is convenient to use a combinatorial description of the crystal tensor product rule which uses certain sequences commonly called signatures.  See, for example, \cite[Sec. 4.4]{Kang}.

For $\lambda \in X(T)$ and $1 \leq j \leq m+n $, define the \emph{$j$-residue of $\lambda$} to be 
\begin{align}\label{E:residuedef}
r_{j}(\lambda)&=(\lambda + \vartheta, \varepsilon_{j})\\
&=\begin{cases} (\lambda + \rho -\varepsilon_{j}, \varepsilon_{j}), & \text{ if $\p{v}_{j}=\0;$}\\
                (\lambda + \rho, \varepsilon_{j}), & \text{ if $\p{v}_{j}=\1.$} \notag
\end{cases} 
\end{align}
Fix $r \in \Z /p\Z$ and $\lambda \in X(T).$ The \emph{$r$-signature} of $\lambda \in X(T)$ is $$\sigma_{r}(\lambda)=(\sigma_{1},  \dotsc , \sigma_{m+n}),$$ where 
\begin{equation}\label{E:sigmadef}
\sigma_{i}=\begin{cases} +, &\text{ if } r_{i}(\lambda+\varepsilon_{i})\equiv r \modp ;\\
                         -, &\text{ if } r_{i}(\lambda) \equiv r \modp ;\\
                         0, &\text{ else.}    
\end{cases}
\end{equation}    Given the $r$-signature of $\lambda$ we form the \emph{reduced $r$-signature} by successively replacing $-+$ pairs with $00$ (where the $-$ and $+$ may be separated by zeros, which are ignored) until no $-$ appears to the left of a $+.$  Given an $r$-signature $\sigma= (\sigma_{1}, \dots , \sigma_{m+n}),$ we write $\tilde{\sigma}=(\tilde{\sigma}_{1}, \dots , \tilde{\sigma}_{m+n})$ for the reduced $r$-signature.  In particular, given $\lambda \in X(T)$ we write $\tilde{\sigma}_{r}(\lambda)$ for the  reduced $r$-signature of $\lambda.$  We then define $\tilde{e}^{*}_{r}, \tilde{f}^{*}_{r}: X(T) \to X(T) \sqcup \{0 \}$ by 
\begin{align*}
\tilde{e}^{*}_{r}(\lambda)&= \begin{cases}   \lambda-\varepsilon_{j}, &\text{ if $1 \leq j \leq m+n$ is the position of the leftmost $-$ in $\tilde{\sigma}_{r}(\lambda);$} \\
                                             0, &\text{ if there are no $-$'s in $\tilde{\sigma}_{r}(\lambda);$} 
\end{cases} \\
\tilde{f}^{*}_{r}(\lambda)&= \begin{cases} \lambda+\varepsilon_{j}, &\text{ if $1 \leq j \leq m+n$ is the position of the rightmost $+$ in $\tilde{\sigma}_{r}(\lambda)$;}\\
                              0, &\text{ if there are no $+$'s in $\tilde{\sigma}_{r}(\lambda)$.} \\
\end{cases}
\end{align*}  
We also define 
\begin{align*}
\varepsilon^{*}_{r}(\lambda) &= \text{max}(a \geq 0: (\tilde{e}_{r}^{*})^{a}(\lambda) \neq 0 )= \text{the total number of $-$'s in $\tilde{\sigma}_{r}(\lambda)$,}\\
\varphi^{*}_{r}(\lambda) &= \text{max}(a \geq 0: (\tilde{f}_{r}^{*})^{a}(\lambda) \neq 0 )= \text{the total number of $+$'s in $\tilde{\sigma}_{r}(\lambda)$.}
\end{align*}  The weight function is given by 
\begin{equation}\label{E:weightdef}
\wt(\lambda) = \sum_{i =1}^{m+n} (-1)^{\p{v}_{i}}\gamma_{(\lambda + \rho, \varepsilon_{i})}.
\end{equation}
  Taken together the datum 
\begin{equation}\label{E:crystalonX(T)}
\left(X(T), \tilde{e}_{r}^{*}, \tilde{f}_{r}^{*}, \varepsilon_{r}^{*}, \varphi_{r}^{*}, \wt \right)
\end{equation}
is the crystal of interest in the present work.  We emphasize that this crystal structure on $X(T)$ depends on (but only on) the sequence of parities $\p{v}_{1}, \dots , \p{v}_{m+n}$ which we fixed at the beginning.

\subsection{Main Results}\label{SS:mainresults}  We can now summarize our main results. Namely, that the crystal on $X(T)$ given in \eqref{E:crystalonX(T)} describes aspects of the category $\Op.$

In section~\ref{S:centralcharacters} we prove that the function $\wt$ given as part of the crystal structure on $X(T)$ partitions the central characters of $\Dist{G}$ arising from the irreducible supermodules of category $\Op$.  In particular if $L(\lambda)$ and $L(\mu)$ have the same central character, then $\wt(\lambda)=\wt(\mu).$  As a consequence we obtain the following linkage principle.

\begin{theorem}\label{T:linkageprinciple}  Let $\lambda, \mu \in X(T).$  If
\[
\Ext^{1}_{\dist{}}(L(\lambda), L(\mu)) \neq 0,
\]
then $\wt(\lambda)=\wt(\mu).$
\end{theorem}

To continue we need to define certain translation functors on category $\Op.$ These functors should be compared with the translation functors defined by Jantzen \cite[II.7]{J1}, Brundan and Kleshchev \cite{BK2}, and Brundan \cite{B0}.  For $\nu \in P,$ define $\Og$ to be the full subcategory of $\Op$ of all modules with all their irreducible subquotients isomorphic to $L(\lambda)$ for some $\lambda \in X(T)$ with $\wt(\lambda)=\nu.$  	Let $M$ be a $\Dist{G}$-supermodule lying in category $\Op.$  Fix $\lambda \in X(T).$ Since $\mathcal{Z}(\Dist{G}),$ the center of $\Dist{G},$ leaves $M_{\lambda}$ invariant, we can view it as a commuting family of endomorphisms of the finite dimensional superspace $M_{\lambda}.$  Consequently, we have the direct sum decomposition as superspaces
\[
M_{\lambda} = \bigoplus M_{\lambda}^{\chi}
\] where the sum runs over all central characters, $\chi: \mathcal{Z}(\Dist{G}) \to k,$ and where 
\[
M_{\lambda}^{\chi} = \{m \in M_{\lambda}: (z-\chi(z))^{N}m=0 \text{ for every $z \in \mathcal{Z}(\Dist{G})$ and $N >>0$} \}. 
\]
Consequently, since $M=\bigoplus_{\lambda \in X(T)}M_{\lambda},$ we have the direct sum decomposition as $\Dist{G}$-supermodules
\[
M = \bigoplus M^{\chi}
\] where again sum runs over all central characters, $\chi: \mathcal{Z}(\Dist{G}) \to k,$ and where
\[
M^{\chi} = \{m \in M: (z-\chi(z))^{N}m=0 \text{ for every $z \in \mathcal{Z}(\Dist{G})$ and $N >>0$ } \}.\]

However, by our above remarks we can rewrite the decomposition as 
\[
M = \bigoplus M^{\nu}
\] where the sum runs over all $\nu \in P$ and $M^{\nu}$ is a $\Dist{G}$-supermodule lying in category $\Og.$  That is, we have 
\begin{equation}\label{E:directsumdecomp}
\Op = \bigoplus_{\nu \in P} \Og.
\end{equation}
Let $\pr_{\nu}$ be the projection functor from $\Op$ to $\Og$.  For $r \in \Z/p\Z$ we now define the functors 
\begin{equation}\label{E:transfunctordefs}  E_{r}, F_{r}: \Op  \to \Op. 
\end{equation}  By additivity, it suffices to define them on objects in $\Og.$  Let $M \in \Og,$ then 
\begin{equation}
\begin{array}{lcl}
E_{r}M:= \pr_{\nu + \gamma_{r}-\gamma_{r+1}} ( M \otimes V^{*} ) &\text{and}& F_{r}M:=\pr_{\nu-\gamma_{r}+\gamma_{r+1}}(M \otimes V),
\end{array}
\end{equation} where $V$ denotes the natural $G$-supermodule and $V^{*}$ denotes its dual.
On a morphism $\varphi: M \to N$, $E_{r}\varphi$ and $F_{r}\varphi$ are the restriction of the map $\varphi \otimes 1.$

 In section~\ref{S:translationfunctors} we prove that the action of these translation functors on irreducible supermodules is regulated by the crystal structure on $X(T).$  Namely, we have the following theorem: 
\begin{theorem}\label{T:transfunctorsonirreps}  Let $\lambda \in X(T)$ and $r \in \Z/p\Z.$
\begin{enumerate}
\item [(i)]  $E_{r}L(\lambda) \neq 0$ if and only if $\varepsilon_{r}^{*}(\lambda) \neq 0,$ in which case it is a self-dual indecomposable module with irreducible socle and cosocle both isomorphic to $L(\mu)$ where $\mu=\tilde{e}^{*}_{r}(\lambda).$  Moreover, $E_{r}L(\lambda)$ is irreducible if and only if $\varepsilon_{r}^{*}(\lambda)=1.$ 
\item [(ii)]   $F_{r}L(\lambda) \neq 0$ if and only if $\varphi_{r}^{*}(\lambda) \neq 0,$ in which case it is a self-dual indecomposable module with irreducible socle and cosocle both isomorphic to $L(\mu)$ where $\mu=\tilde{f}^{*}_{r}(\lambda).$  Moreover, $F_{r}L(\lambda)$ is irreducible if and only if $\varphi_{r}^{*}(\lambda)=1.$ 
\end{enumerate}
\end{theorem}
As a corollary to the previous theorem we also obtain an explicit description of the socle of $L(\lambda)\otimes V^{*}$ and $L(\lambda) \otimes V$ and combinatorial criterion for when these supermodules are semisimple.

\section{Central Characters}\label{S:centralcharacters}

\subsection{Some Central Elements}\label{S:prelimdefs}

Following Sergeev \cite{sergeev2} (where the characteristic zero case was considered), we define certain central elements of $\Dist{G}$.  Recall from subsection~\ref{S:algebraofdistributions} that $\{ e_{i,j}: 1 \leq i,j \leq m+n\}$ is the usual homogeneous basis for the Lie superalgebra associated to $G.$   Define $x_{k,l}^{[r]} \in \dist{}$ inductively as follows:
\begin{align}\label{E:xdef}
x_{k,l}^{[1]}&=e_{k,l} \\
x_{k,l}^{[r]}&=\sum_{s =1}^{m+n}(-1)^{\p{v}_{s}}e_{k,s}x_{s,l}^{[r-1]}, \mbox{ for } r > 1.
\end{align}

\begin{lemma}\label{L:basicfacts}Let $1 \leq  i,j,k,l \leq m+n$ and $r \in \Z, r \geq 1.$  Then,
\begin{enumerate}
\item [(i)] The $x_{k,l}^{[r]}$ are homogeneous of degree $\p{v}_{k}+\p{v}_{l}$, 
\item [(ii)] $[e_{i,j},x_{k,l}^{[r]}]=\delta_{j,k}x_{i,l}^{[r]}-(-1)^{(\p{v}_{i}+\p{v}_{j})(\p{v}_{k}+\p{v}_{l})}\delta_{i,l}x_{k,j}^{[r]}$,
\item [(iii)] The elements $\tilde{Z}_{r}:=\sum_{k =1}^{m+n} x_{k,k}^{[r]}$ for $r \geq 1$ are central.
\end{enumerate}
\end{lemma}

\begin{proof}  Each of the statements is a straightforward induction on $r.$
\end{proof}

For $r \in \Z$ let
\begin{equation}\label{E:zrdef}
Z_{r}=\tilde{Z}_{r}-(-1)^{r}\sum (-1)^{\p{v}_{k_{1}}+\dots +\p{v}_{k_{r}}}\vartheta_{k_{r}},
\end{equation} where the sum runs over all $1 \leq k_{1}< \dots < k_{r} \leq m+n$ and $\vartheta_{i}$ is as in \eqref{E:thetadef}.  Since $Z_{r}$ differs from $\tilde{Z}_{r}$ by a scalar, these elements are still central and generate the same subalgebra of $\Dist{G}.$  We remark that if $\operatorname{char} k=0$ one can use the results of \cite{sergeev3} to prove that these elements in fact \emph{generate} $\mathcal{Z}(\Dist{G}).$  

\subsection{The Linkage Principle}\label{SS:linkageprinciple}

Multiplication by an element of $\mathcal{Z}(\Dist{G})$ defines an endomorphism of $M(\lambda)$ ($\lambda \in X(T)$) so takes the canonical generator $v_{\lambda} \in M(\lambda)_{\lambda}$ to scalar multiple of itself.  From this we conclude that the elements of $\mathcal{Z}(\Dist{G})_{\0}$ must act by scalars on $M(\lambda)$ and any of its subquotients.  We now prove that the even central element $Z_{r}$ defined in the previous subsection acts on $M(\lambda)$ by the scalar $Z_{r}(\lambda)$ given below.  

For $1 \leq i \leq m+n$ define 
\begin{equation}\label{E:ydefandylambdadef}
r_{i}=(-1)^{\p{v}_{i}}(e_{i,i}+\vartheta_{i}),
\end{equation} where $\vartheta_{i}$ is as defined in \eqref{E:thetadef}.
Observe that if $M$ is a $\dist{}$-supermodule and $v \in M_{\lambda}$, then $r_{i}v=r_{i}(\lambda)v,$ the $i$th residue of $\lambda$ \eqref{E:residuedef}.  Given $\lambda \in X(T)$ define the integer 
\begin{equation}\label{E:zrlambdadef}
Z_{r}(\lambda) = \sum_{s=1}^{r}(-1)^{s-1}\sum_{}(-1)^{\p{v}_{k_{1}}+\dotsb + \p{v}_{k_{s}}}r_{k_{1}}(\lambda)^{a_{1}}\dotsb r_{k_{s}}(\lambda)^{a_{s}}
\end{equation}
where the unmarked sum runs over all $1 \leq k_{1}< \cdots < k_{s} \leq m+n$ and nonnegative integers $a_{1}, \dotsc , a_{s}$ such that $a_{1}+\dotsb +a_{s}=r-s+1$.

\begin{lemma}\label{T:scalar}  If $M$ is a $\dist{}$-supermodule and $v \in M_{\lambda}$  (for some $\lambda \in X(T)$) is annihilated by all $e_{i,j}$ when $1 \leq i< j \leq m+n$, then $Z_{r}v=Z_{r}(\lambda)v$.  
\end{lemma}

\begin{proof}  Given our assumption that $v$ is annihilated by all $e_{i,j}$ when $i < j$, it suffices to work work modulo the left ideal generated by these elements.  We shall write $\equiv$ for congruence modulo this ideal.  We shall prove the statement via several intermediate claims:

\vspace{1mm}
\noindent{\bf Claim 1.} {\em $x^{[r]}_{k,l} \equiv 0$ for each $1 \leq k <l \leq m+n$, and $r \geq 1$.}
\vspace{1mm}

We prove this by inducting on $r$ with the base case being clear.  For $r > 1$ by the induction hypothesis, Lemma~\ref{L:basicfacts}(ii), and \eqref{E:xdef}  we have:

\begin{align*} x^{[r]}_{k,l} & = \sum_{s=1}^{m+n}(-1)^{\p{v}_{s}}e_{k,s}x_{s,l}^{[r-1]} \\
& \equiv \sum_{s \geq l} (-1)^{\p{v}_{s}}e_{k,s}x_{s,l}^{[r-1]} \\
& = \sum_{s \geq l} (-1)^{\p{v}_{s}}\left(x_{k,l}^{[r-1]}+(-1)^{(\p{v}_{k}+\p{v}_{s})(\p{v}_{s}+\p{v}_{l})}x_{s,l}^{[r-1]}e_{k,s} \right) \\
& \equiv 0
\end{align*}

\vspace{1mm}
\noindent{\bf Claim 2:} {\em 
\[
x_{k,k}^{[r]} \equiv r_{k}x_{k,k}^{[r-1]}- (-1)^{\p{v}_{k}} \sum_{s > k } x_{s,s}^{[r-1]}.
\]
for all $1 \leq k \leq m+n$ and $r \geq 1$ (where we define $x_{s,s}^{[0]}=(-1)^{\p{v}_{s}}$).}
\vspace{1mm}

We do a direct calculation using (\ref{E:xdef}), Claim 1, and Lemma \ref{L:basicfacts} (ii):

\begin{align*} x_{k,k}^{[r]} &= \sum_{s =1}^{m+n} (-1)^{\p{v}_{s}}e_{k,s}x_{s,k}^{[r-1]} \\
&\equiv (-1)^{\p{v}_{k}}e_{k,k}x_{k,k}^{[r-1]}+\sum_{s>k}(-1)^{\p{v}_{s}}e_{k,s}x_{s,k}^{[r-1]} \\
& \equiv (-1)^{\p{v}_{k}}h_{k}x_{k,k}^{[r-1]}  +\sum_{s>k}(-1)^{\p{v}_{s}}\left(x_{k,k}^{[r-1]} - (-1)^{\p{v}_{k}+\p{v}_{s}}x_{s,s}^{[r-1]}+(-1)^{\p{v}_{k} + \p{v}_{s}}x_{s,k}^{[r-1]}e_{k,s} \right)\\
&\equiv (-1)^{\p{v}_{k}}h_{k}x_{k,k}^{[r-1]} + \sum_{s>k}(-1)^{\p{v}_{s}}x_{k,k}^{[r-1]}-\sum_{s>k}(-1)^{\p{v}_{k}}x_{s,s}^{[r-1]} \\
&=r_{k}x_{k,k}^{[r-1]}-(-1)^{\p{v}_{k}}\sum_{s>k}x_{s,s}^{[r-1]}.
\end{align*}

\vspace{1mm}
\noindent{\bf Claim 3:} {\em
\begin{equation}\label{E:nastysum}
x_{k,k}^{[r]}\equiv \sum_{s=1}^{r}(-1)^{s-1}\sum (-1)^{\p{v}_{k_{1}} + \dotsb + \p{v}_{k_{s-1}}}r_{k_{1}}^{a_{1}}\dotsb r_{k_{s}}^{a_{s}}h_{k_{s}},
\end{equation} for $r \geq 1$;
where the second sum is over all $k = k_{1} < k_{2} < \dotsb < k_{s}$ and $a_{1}, \dotsc , a_{s} \in \mathbb{Z}_{\geq 0}$ such that $a_{1}+ \dotsb +a_{s}=r-s.$  }
\vspace{1mm}

We induct on $r \geq 1$ with the case $r=1$ being clear.  Let $r >1$, then by Claim 2 and the induction hypothesis we have

\begin{align*} x_{k,k}^{[r]} &\equiv  r_{k}x_{k,k}^{[r-1]}- (-1)^{\p{v}_{k}}\sum_{s> k} x_{s,s}^{[r-1]} \\
 &\equiv  r_{k}\left(\sum_{t=1}^{r-1}(-1)^{t-1}\sum (-1)^{\p{v}_{k_{1}} + \dotsb + \p{v}_{k_{t-1}}}r_{k_{1}}^{a_{1}}\dotsb r_{k_{t}}^{a_{t}}h_{k_{t}} \right) \\
&  \hspace{1in} -(-1)^{\p{v}_{k}}\sum_{s>k}\left(\sum_{u=1}^{r-1}(-1)^{u-1}\sum (-1)^{\p{v}_{k_{1}} + \dotsb + \p{v}_{k_{u-1}}}r_{k_{1}}^{a_{1}}\dotsb r_{k_{u}}^{a_{u}}h_{k_{u}} \right) \\
&\equiv \sum_{t=1}^{r-1}(-1)^{t-1}\sum (-1)^{\p{v}_{k_{1}} + \dotsb + \p{v}_{k_{t-1}}}r_{k_{1}}^{a_{1}+1}\dotsb r_{k_{t}}^{a_{t}}h_{k_{t}}\\
& \hspace{1in} -(-1)^{\p{v}_{k}}\sum_{s>k}\left(\sum_{u=1}^{r-1}(-1)^{u-1}\sum (-1)^{\p{v}_{k_{1}} + \dotsb + \p{v}_{k_{u-1}}}r_{k_{1}}^{a_{1}}\dotsb r_{k_{u}}^{a_{u}}h_{k_{u}} \right)
\end{align*}
Which, one observes, is equal to the double sum given in \eqref{E:nastysum}.

Consequently, we have 
\[
\tilde{Z}_{r} \equiv \sum_{s=1}^{r}(-1)^{s-1}\sum (-1)^{\p{v}_{k_{1}} + \dotsb + \p{v}_{k_{s-1}}}r_{k_{1}}^{a_{1}}\dotsb r_{k_{s}}^{a_{s}}h_{k_{s}}
\] for $r \geq 1$;
where the second sum is over all $k_{1} < k_{2} < \dotsb < k_{s}$ and $a_{1}, \dotsc , a_{s} \in \mathbb{Z}_{\geq 0}$ such that $a_{1}+ \dotsb +a_{s}=r-s$.

\vspace{1mm}
\noindent{\bf Claim 4:} {\em  
\begin{equation}\label{E:zrformula}
Z_{r} \equiv \sum_{s=1}^{r}(-1)^{s-1}\sum (-1)^{\p{v}_{k_{1}}+\dotsb +\p{v}_{k_{s}}}r_{k_{1}}^{a_{1}}\dotsb r_{k_{s}}^{a_{s}},
\end{equation}
where the unmarked sum runs over all $k_{1}<\dotsb <k_{s}$ and nonnegative integers $a_{1}, \dotsc , a_{s}$ such that $a_{1}+\dotsb + a_{s}=r-s+1$.}
\vspace{1mm}

Using (\ref{E:zrdef}) and Claim 2 we have

\begin{align*} Z_{r} &= \tilde{Z}_{r}-(-1)^{r}\sum_{k_{1} < \dotsb < k_{r}}(-1)^{\p{v}_{k_{1}}+\dotsb +\p{v}_{k_{r}}}\vartheta_{k_{r}}\\
&\equiv \sum_{s=1}^{r}(-1)^{s-1}\sum (-1)^{\p{v}_{k_{1}} + \dotsb + \p{v}_{k_{s-1}}}r_{k_{1}}^{a_{1}}\dotsb r_{k_{s}}^{a_{s}}h_{k_{s}} -(-1)^{r}\sum_{k_{1} < \dotsb < k_{r}}(-1)^{\p{v}_{k_{1}}+\dotsb +\p{v}_{k_{r}}}\vartheta_{k_{r}}
\intertext{however $r_{k_{s}}=(-1)^{\p{k_{s}}}(h_{k_{s}}+\vartheta_{k_{s}})$ so a substitution yields:}
&= \sum_{s=1}^{r}(-1)^{s-1}\sum (-1)^{\p{v}_{k_{1}} + \dotsb + \p{v}_{k_{s}}}r_{k_{1}}^{a_{1}}\dotsb r_{k_{s}}^{a_{s}+1}\\
& \hspace{1in}- \sum_{s=1}^{r}(-1)^{s-1}\sum (-1)^{\p{v}_{k_{1}} + \dotsb + \p{v}_{k_{s-1}}}r_{k_{1}}^{a_{1}}\dotsb r_{k_{s}}^{a_{s}}\vartheta_{k_{s}}\\
& \hspace{2.5in}-(-1)^{r}\sum_{k_{1} < \dotsb < k_{r}}(-1)^{\p{v}_{k_{1}}+\dotsb +\p{v}_{k_{r}}}\vartheta_{k_{r}}\\
&= \sum_{s=1}^{r}(-1)^{s-1}\sum (-1)^{\p{v}_{k_{1}} + \dotsb + \p{v}_{k_{s}}}r_{k_{1}}^{a_{1}}\dotsb r_{k_{s}}^{a_{s}+1}\\
& \hspace{1in}+ \sum_{s=1}^{r}(-1)^{s}\sum (-1)^{\p{v}_{k_{1}} + \dotsb + \p{v}_{k_{s}}}r_{k_{1}}^{a_{1}}\dotsb r_{k_{s}}^{a_{s}}\left(\sum_{t>k_{s}}(-1)^{\p{v}_{t}}r_{t}^{0} \right)\\
&\hspace{2.5in}-(-1)^{r}\sum_{k_{1} < \dotsb < k_{r}}(-1)^{\p{v}_{k_{1}}+\dotsb +\p{v}_{k_{r}}}\\
&\equiv \sum_{s=1}^{r}(-1)^{s-1}\sum (-1)^{\p{v}_{k_{1}} + \dotsb + \p{v}_{k_{s}}}r_{k_{1}}^{a_{1}}\dotsb r_{k_{s}}^{a_{s}+1}\\
& \hspace{1in} + \sum_{s=1}^{r}(-1)^{s}\sum_{\substack{k_{1}<\dotsb <k_{s}<t \\ a_{1}+ \dotsb +a_{s}=r-s }} (-1)^{\p{v}_{k_{1}} + \dotsb + \p{v}_{k_{s}}+\p{v}_{t}}r_{r_{1}}^{a_{1}}\dotsb r_{k_{s}}^{a_{s}}r_{t}^{0} \\
&\hspace{2.5in}-(-1)^{r}\sum_{k_{1} < \dotsb < k_{r}}(-1)^{\p{v}_{k_{1}}+\dotsb +\p{v}_{k_{r}}}
\end{align*}
We now observe that the first sum is \eqref{E:zrformula} in all cases when $a_{s}\geq 1$, the second sum (after reindexing) is \eqref{E:zrformula} for $a_{s}=0$ and $s=2,\dotsc, r$ plus one additional term which exactly cancels the third sum.  It remains to note that the final possible case (when $s=1$ and $a_{s}=0$) in fact does not occur in \eqref{E:zrformula}. 

Finally, since $r_{t}v=r_{t}(\lambda)v$, Claim 4 exactly implies the statement given in the lemma.
\end{proof}

For $\lambda \in X(T)$, define $G_{\lambda}(t) \in k[[t^{-1}]]$ by 
\begin{equation}\label{E:gdef}
G_{\lambda}(t) = 1-\sum_{r \geq 1}Z_{r}(\lambda)t^{-(r+1)},
\end{equation}
and for each $r \in \Z/p\Z$ define 
\begin{align}\label{E:ABabdef}
A_{r}(\lambda) &= \left| \setof{i}{r_{i}(\lambda +\varepsilon_{i}) \equiv r \: (\mbox{mod } p)} \right|,\\
B_{r}(\lambda) &= \left| \setof{i}{r_{i}(\lambda) \equiv r \: (\mbox{mod } p)} \right|, \notag 
\end{align}

\begin{lemma}\label{L:linkagecombinatorial}  Let $\lambda, \mu \in X(T).$  If $\lambda=\sum_{i =1}^{m+n}\lambda_{i}\varepsilon_{i},$ then let $l(\lambda)$ denote $\sum_{i =1}^{m+n} \lambda_{i}.$   The following are equivalent:

\begin{enumerate}
\item [(i)] $l(\lambda)=l(\mu)$ and $Z_{s}(\mu) \equiv Z_{s}(\mu) \modp$ for all $s \in \Z_{\geq 0}$;
\item [(ii)]  $l(\lambda)=l(\mu)$ and $G_{\lambda}(t)=G_{\mu}(t)$;
\item [(iii)]  $l(\lambda)=l(\mu)$ and $A_{r}(\lambda) - B_{r}(\lambda) = A_{r}(\mu) - B_{r}(\mu)$ for all $r\in \Z/p\Z$;
\item [(iv)] $\wt(\lambda)=\wt(\mu);$
\end{enumerate} where $\wt: X(T) \to P$ is the weight function defined in \eqref{E:weightdef}.
\end{lemma} 

\begin{proof}  $(i) \Leftrightarrow (ii)$ is immediate from \eqref{E:gdef}.

To prove $(ii) \Leftrightarrow (iii)$, we observe that 
\begin{align*}
G_{\lambda}(t) &= \prod_{i =1}^{m+n}\left(1-\sum_{r\geq 1}(-1)^{\p{v}_{i}}r_{i}(\lambda)^{r-1}t^{-r} \right) \\
& = \prod_{i=1}^{m+n}\left(1-(-1)^{\p{v}_{i}}t^{-1}\sum_{r \geq 1}(r_{i}(\lambda)t^{-1})^{r-1} \right) \\
& = \prod_{i=1}^{m+n}\left(1-(-1)^{\p{v}_{i}}t^{-1}\left(\frac{1}{1-r_{i}(\lambda)t^{-1}} \right) \right) \\
& = \prod_{i=1}^{m+n}\left(\frac{1-r_{i}(\lambda)t^{-1}-(-1)^{\p{v}_{i}}t^{-1}}{1-r_{i}(\lambda)t^{-1}} \right) \\
&= \prod_{i=1}^{m+n}\left(\frac{t-r_{i}(\lambda + \varepsilon_{i})}{t-r_{i}(\lambda)} \right). 
\end{align*}
Comparing the multiplicity of zeros and poles we see $G_{\lambda}(t)=G_{\mu}(t)$ if and only if $A_{r}(\lambda) - B_{r}(\lambda) = A_{r}(\mu) - B_{r}(\mu)$ for $r\in \Z/p\Z.$

To prove $(iii) \Leftrightarrow (iv)$, we need to analyze the weight function more closely.  In characteristic zero, we have 
\begin{align*}
\langle (-1)^{\p{v}_{i}}\gamma_{(\lambda+\rho, \varepsilon_{i})}, \alpha_{r}\rangle &=\langle (-1)^{\p{v}_{i}}\gamma_{(\lambda+\rho, \varepsilon_{i})}, \gamma_{r}-\gamma_{r+1}\rangle\\
&= \begin{cases} 1, &\text{ if } r_{i}(\lambda+\varepsilon_{i})=r;\\
                 -1, &\text{ if } r_{i}(\lambda)=r;\\
                 0, &\text{ otherwise.}  
\end{cases}
\end{align*}
That is, $\langle \wt(\lambda),\alpha_{r} \rangle = A_{r}(\lambda)-B_{r}(\lambda)$ for all $r \in \Z.$  Now if $\lambda, \mu \in X(T)$ with $\wt(\lambda)=\wt(\mu),$ then 
\[
0=\sum_{r \in \Z}\langle \wt(\lambda)-\wt(\mu), r\gamma_{r} \rangle=l(\lambda)-l(\mu),
\] and, for all $r\in \Z,$
\[
0=\langle \wt(\lambda)-\wt(\mu), \alpha_{r} \rangle = (A_{r}(\lambda)-B_{r}(\lambda)) -(A_{r}(\mu)-B_{r}(\mu)).
\]  That is, $(iv) \Rightarrow (iii).$  To prove the converse, first observe that $\wt(\lambda)-\wt(\mu)$ lies in the $\Z$-span of the $\alpha_{r}$'s so we can fix $N>0$ so that $\wt(\lambda)-\wt(\mu)$ lies in the $\Z$-span of $\{\alpha_{-N}, \dots ,\alpha_{N} \}.$  If $(iii)$ holds, then $\langle \wt(\lambda)-\wt(\mu), \alpha_{r}\rangle=0$ for $-N \leq r \leq N.$  That is, $\wt(\lambda)-\wt(\mu)$ lies in the radical of the form restricted to this sublattice.  However, the form is nondegenerate on this sublattice.  This implies $\wt(\lambda)=\wt(\mu).$

Now consider the case when $\operatorname{char} k > 0$.  For $1 \leq i \leq m+n,$ write $(\lambda + \rho, \varepsilon_{i})=pb_{i}+r_{i}$ where $b_{i} \in \Z,$ and $r_{i}=1, \dots , p.$  Then, 
\begin{align*}
\wt(\lambda)&=\sum_{i =1}^{m+n} (-1)^{\p{v}_{i}}\gamma_{(\lambda+\rho, \varepsilon_{i})}\\
&=\sum_{i =1}^{m+n} (-1)^{\p{v}_{i}}(\gamma_{r_{i}} - b_{i}\delta)\\
&=\sum_{i=1}^{m+n}(-1)^{\p{v}_{i}}(\Lambda_{r_{i}}-\Lambda_{r_{i}-1}) - \sum_{i =1}^{m+n}(-1)^{\p{v}_{i}}b_{i}\delta\\
&=\sum_{i=1}^{m+n}(\Lambda_{r_{i}(\lambda+\varepsilon_{i})}-\Lambda_{r_{i}(\lambda)}) - \sum_{i =1 }^{m+n}(-1)^{\p{v}_{i}}b_{i}\delta\\
&=\sum_{r \in \Z/p\Z}(A_{r}(\lambda)-B_{r}(\lambda))\Lambda_{r} - \sum_{i =1}^{m+n}(-1)^{\p{v}_{i}}b_{i}\delta.
\end{align*} In particular, we see that $\langle \wt(\lambda), \alpha_{r} \rangle = A_{r}(\lambda)-B_{r}(\lambda)$ for $r \in \Z/p\Z.$   Let $K=-\sum_{r \in \Z/p\Z}\Lambda_{r} \in P.$  Then $\langle \gamma_{a}, K \rangle =a$ for all $a \in \Z.$  Consequently, 
\[
\langle \wt(\lambda), K \rangle=\sum_{i =1}^{m+n}(\lambda_{i}+\rho_{i})=l(\lambda)+\sum_{i =1}^{m+n}\rho_{i}.
\]  Let $\lambda, \mu \in X(T).$  From the above remarks, $(iv)$ holds if and only if $\langle \wt(\lambda)-\wt(\mu), \alpha_{r}\rangle=0$ for $r\in \Z/p\Z$ and $\langle \wt(\lambda)-\wt(\mu), K \rangle= 0.$  However the $\alpha_{r}$'s and $K$ together span $\mathfrak{h}$ and the form is nondegenerate on $\mathfrak{h}.$  Consequently, $(iv)$ holds if and only if $\wt(\lambda)=\wt(\mu).$ 

\end{proof}

We can now deduce the main results of this section. 

\begin{theorem}\label{T:linkageprinciple1}  Let $\lambda, \mu\in X(T).$   If $L(\mu)$ is a subquotient of $M(\lambda),$ then $\wt(\lambda)=\wt(\mu).$
\end{theorem} 
\begin{proof}  If $L(\mu)$ is a composition factor of $M(\lambda),$ then $\mathcal{Z}(\Dist{G})_{\0}$ acts by the same scalars on both.  If we set
\[
d_{r}=\binom{\sum_{i =1}^{m+n} e_{i,i}}{r}
\] for $r \in \Z_{\geq 0},$ a direct calculation verifies that these even elements of $\Dist{G}$ are central.  Furthermore, $d_{l(\lambda)}$ acts by $1$ and $d_{r}$ acts by $0$ for all $r > l(\lambda)$ (c.f. \cite[Remark 2, Sec. 3.8]{Carter}).  Consequently, $l(\lambda)=l(\mu).$ Additionally, we have $Z_{r}(\lambda)=Z_{r}(\mu)$ for all $r \in \Z.$  Taken together with Lemma~\ref{L:linkagecombinatorial} this implies $\wt(\lambda)=\wt(\mu).$
\end{proof}

\begin{lemma}\label{L:exttriviality}  Let $\lambda\in X(T)$ and $D$ be a $\Dist{G}$-supermodule such that if $\mu\in X(T)$ with $\lambda < \mu$ in the dominance order, then $D_{\mu}=0$.  Then
\[
\Ext_{\Dist{G}}^{1}(M(\lambda), D) =0.
\]  Moreover, for arbitrary $\lambda, \mu\in X(T)$ we have 
\[
\Ext_{\Dist{G}}^{1}(M(\lambda), W(\mu)) =0.
\]
\end{lemma}

\begin{proof}  Consider a short exact sequence 
\[
0 \to D \to N \to M(\lambda) \to 0.
\] Let $v'_{\lambda} \in N$ be a homogeneous preimage of $v_{\lambda},$ the canonical generator of $M(\lambda).$  By the weight assumption we have that $v'_{\lambda}$ is a primitive vector of weight $\lambda.$  By the universal property of $M(\lambda)$ there is a homomorphism which maps $v_{\lambda} \mapsto v'_{\lambda},$ providing a splitting of the sequence.  Hence $\Ext^{1}_{\Dist{G}}(M(\lambda), D) = 0.$

Now if $\lambda, \mu \in X(T)$ are arbitrary, by contravariant duality we have 
\[
\Ext_{\Dist{G}}^{1}(M(\lambda), W(\mu))  \cong \Ext_{\Dist{G}}^{1}(M(\mu), W(\lambda)),
\] so we can assume without loss that  $\lambda \not< \mu$.  The result then immediately follows. 
\end{proof}

\begin{lemma}\label{L:radicalhomandext} If $\lambda, \mu\in X(T)$ with $\lambda \not< \mu,$ then 
\[
\Hom_{\Dist{G}}(\Rad M(\lambda), L(\mu)) \cong \Ext_{\Dist{G}}^{1}(L(\lambda), L(\mu)).
\]  In particular, for all $\lambda \in X(T)$ 
\[
\Ext_{\Dist{G}}^{1}(L(\lambda), L(\lambda))=0.
\]
\end{lemma}

\begin{proof} Let $R$ denote $\Rad M(\lambda).$  Consider the short exact sequence 
\[
0 \to R \to M(\lambda) \to L(\lambda) \to 0.
\]  Applying the functor $\Hom_{\Dist{G}}(-, L(\mu))$ we obtain the long exact sequence 
\begin{multline*}
0 \to \Hom_{\Dist{G}}(L(\lambda), L(\mu)) \to \Hom_{\Dist{G}}(M(\lambda), L(\mu)) \to\\
 \to\Hom_{\Dist{G}}(R, L(\mu))
\to \Ext_{\Dist{G}}^{1}(L(\lambda), L(\mu)) \to \\
\to \Ext_{\Dist{G}}^{1}(M(\lambda), L(\mu)) \to \dots .
\end{multline*}
We note that the map $\Hom_{\Dist{G}}(L(\lambda), L(\mu)) \to \Hom_{\Dist{G}}(M(\lambda), L(\mu))$ is an isomorphism. Since $\lambda \not< \mu,$ by Lemma~\ref{L:exttriviality} we have that $\Ext_{\Dist{G}}^{1}(M(\lambda), L(\mu))=0.$  The result immediately follows.  The triviality of $\Ext_{\Dist{G}}^{1}(L(\lambda),L(\lambda))$ follows by a weight argument.
\end{proof}

\begin{theorem}\label{T:linkageprincipleagain}  Let $\lambda, \mu \in X(T).$  If
\[
\Ext^{1}_{\dist{}}(L(\lambda), L(\mu)) \neq 0,
\]
then $\wt(\lambda)=\wt(\mu).$
\end{theorem}

\begin{proof} By contravariant duality we have
\[
\Ext^{1}_{\Dist{G}}(L(\lambda), L(\mu)) \cong \Ext^{1}_{\Dist{G}}(L(\mu), L(\lambda)),
\] so we may assume without loss that $ \lambda \not< \mu $ in the dominance order.  By Lemma~\ref{L:radicalhomandext}, we can then deduce that $L(\mu)$ is a quotient of $\Rad M(\lambda),$ hence a subquotient of $M(\lambda).$  Finally, one applies Theorem \ref{T:linkageprinciple1}.
\end{proof}

\section{Lowering operators}\label{S:loweringoperators}

In a series of papers (\cite{K2}, \cite{K3}, \cite{K4}) Kleshchev gave certain elements of the algebra of distributions of $SL(n)$ called \emph{lowering operators} which allowed him to prove modular branching rules for $SL(n),$ and via Schur functor arguments, the symmetric group.  Brundan generalized Kleshchev's construction to the case of quantum $GL(n)$.  In this section we develop the appropriate theory for the supergroup $G=GL(V)$ and use them to prove a representation theoretic interpretation of the crystal operators $\tilde{e}_{r}^{*}$ and $\tilde{f}_{r}^{*}.$  

\subsection{A Combinatorial Interlude}\label{S:combinatorialdefinitions}  We now introduce the super analogue of Kleshchev's combinatorial notions of normal, conormal, good, and cogood.  We first define them in terms of the crystal structure on $X(T)$ and then relate them to the appropriate generalization of Kleshchev's original definitions.  Let $r \in \Z/p\Z$ and $\lambda \in X(T).$  

\begin{definition}\label{D:normal}  We define $1 \leq i \leq m+n$ to be \emph{$r$-normal for $\lambda$} if one of the following equivalent conditions hold: 
\begin{enumerate}
\item $(\tilde{e}_{r}^{*})^{a+1}(\lambda)=(\tilde{e}_{r}^{*})^{a}(\lambda) - \varepsilon_{i}$ for some $a \in \Z_{\geq 0};$
\item $i$ is the position of a $-$ in the reduced $r$-signature of $\lambda.$
\end{enumerate}
\end{definition}

\begin{definition}\label{D:good}  We define $1 \leq i \leq m+n$ to be \emph{$r$-good for $\lambda$} if one of the following equivalent conditions hold: 
\begin{enumerate}
\item $\tilde{e}_{r}^{*}(\lambda)=\lambda - \varepsilon_{i};$
\item $i$ is the position of the leftmost $-$ in the reduced $r$-signature of $\lambda.$
\end{enumerate}
\end{definition}

We also have the analogous definitions for $\tilde{f}_{r}^{*}$: 

\begin{definition}\label{D:conormal}  We define $1 \leq i \leq m+n$ to be \emph{$r$-conormal for $\lambda$} if one of the following equivalent conditions hold: 
\begin{enumerate}
\item $(\tilde{f}_{r}^{*})^{a+1}(\lambda)=(\tilde{f}_{r}^{*})^{a}(\lambda) + \varepsilon_{i}$ for some $a \in \Z_{\geq 0};$
\item $i$ is the position of a $+$ in the reduced $r$-signature of $\lambda.$
\end{enumerate}
\end{definition}

\begin{definition}\label{D:cogood}  We define $1 \leq i \leq m+n$ to be \emph{$r$-cogood for $\lambda$} if one of the following equivalent conditions hold: 
\begin{enumerate}
\item $\tilde{f}_{r}^{*}(\lambda)=\lambda + \varepsilon_{i};$
\item $i$ is the position of the rightmost $+$ in the reduced $r$-signature of $\lambda.$
\end{enumerate}
\end{definition}

We call $1 \leq i \leq m+n $ \emph{normal} (resp. \emph{good}, \emph{conormal}, \emph{cogood}) for $\lambda$ if $i$ is $r$-normal (resp. $r$-good, $r$-conormal, $r$-cogood) for $\lambda$ for some $r \in \Z/p\Z.$  

\begin{example}\textnormal{Assume our fixed homogeneous basis for $V$ satisfies $\p{v}_{i}=\1 $ if $1 \leq i \leq n$ and $\p{v}_{i}=\0$ if $n+1 \leq i \leq m+n.$  Consider the crystal structure $(X(T), \tilde{e}_{r}^{*}, \tilde{f}_{r}^{*}, \varepsilon_{r}^{*}, \varphi_{r}^{*}, \wt ).$  Let $p=3,$ $m=3,n=2,$ and 
\[
\lambda=\varepsilon_{1}-\varepsilon_{2}+\varepsilon_{3}+7\varepsilon_{4}+5\varepsilon_{5} \in X(T).
\]
The signatures of $\lambda$ are: $\sigma_{0}(\lambda)=(+,+,-,+,+),$ $\sigma_{1}(\lambda)=(-,-,+,0,0),$ and $\sigma_{2}(\lambda)=(0,0,0,-,-).$  The reduced signatures of $\lambda$ are: $\tilde{\sigma}_{0}(\lambda)=(+,+,0,0,+),$ $\tilde{\sigma}_{1}(\lambda)=(-,0,0,0,0),$ and $\tilde{\sigma}_{2}(\lambda)=(0,0,0,-,-).$  Then for $\lambda:$ $1$ is $1$-normal and $1$-good; $4$ and $5$ are $2$-normal with $4$ being 2-good; $1,2$ and $5$ are $0$-conormal with $5$ being $0$-cogood.}
\end{example} 

We now relate these crystal theoretic definitions to the super analogue of Kleshchev's original definitions.  To do so we require additional notation.  Given $1 \leq i < j\leq m+n,$ let $(i..j)=\{i+1, \dots , j-1 \},$ $(i..j]=\{i+1, \dots , j \},$ $[i..j)=\{i, \dots , j-1 \},$ and $[i..j]=\{i,\dots ,j \}.$  Given $A \subseteq \{1, \dotsc , m+n \}$ let $A_{i..j}=A \cap (i..j).$  Given $A, B \subseteq \{1, \dotsc , m+n \},$ we write $A \downarrow B$ if there exists an injective map $\theta: A \to B$ with $\theta(a) \leq a$ for all $a \in A.$  An equivalent definition is to say $A \downarrow B$ if $|A \cap [1..k]|\leq |B \cap [1..k]|$ for all $1 \leq k \leq m+n.$ This defines a partial order on the subsets of $\{1, \dotsc , m+n \}.$

For $1 \leq i,j \leq m+n,$ $1 \leq k \leq m+n-1,$ and $\lambda \in X(T)$, let 
\begin{align}\label{E:candbdefs}
c_{i,j}(\lambda) &= (\lambda + \vartheta, \varepsilon_{i}-\varepsilon_{j}),\\
b_{i,k}(\lambda) &=(\lambda + \vartheta + \varepsilon_{k+1} , \varepsilon_{i}-\varepsilon_{k+1}), 
\end{align}   Let $1 \leq i <j \leq m+n.$  Define the following subsets of $\{1, \dotsc , m+n \}$:
\begin{align}\label{E:CandBdefs}
C_{i,j}(\lambda) &=\{h \in (i..j) : c_{i,h}(\lambda) \equiv 0  \modp\},\\
B_{i,j}(\lambda) &=\{h \in [i..j) : b_{i,h}(\lambda) \equiv 0  \modp \}. 
\end{align}

Let $\lambda \in X(T)$ and let $1 \leq i, h \leq m+n.$  The following remarks are immediate from the definitions.
\begin{remark} \label{R:hinC} Let $\lambda \in X(T)$ and let $1 \leq i, h \leq m+n,$ then
\begin{enumerate}
\item $h \in C_{i,m+n}(\lambda)$ if and only if $i < h < n$ and $r_{i}(\lambda) \equiv r_{h}(\lambda) \modp;$ 
\item $h \in B_{i,m+n}(\lambda)$ if and only if $i \leq h < n$ and $r_{i}(\lambda) \equiv r_{h+1}(\lambda + \varepsilon_{h+1}) \modp.$ 
\end{enumerate}
\end{remark}

Using the above remarks and the definitions of normal and good, a straightforward combinatorial argument proves the following result.

\begin{lemma} \label{L:crystalnormalcriteria}  Let $\lambda \in X(T)$ and $1 \leq i \leq m+n.$  Then $i$ is normal for $\lambda$ if and only if $B_{i,m+n}(\lambda) \downarrow C_{i,m+n}(\lambda)$.  Furthermore, $i$ is good for $\lambda$ if and only if $i$ is normal and there is no $j$ which is normal, $j <i,$ and $c_{j,i}(\lambda) \equiv 0 \modp.$  
\end{lemma}
We remark that in the purely even case these are the definitions of normal and good used by Brundan.  They are transpose to those given by Kleshchev (see \cite{B2}).

We end this subsection with two combinatorial results which will be useful later.  A straightforward combinatorial argument proves the following lemma.

\begin{lemma}\label{L:normalplusconormalequalsgood}  Given $1 \leq i \leq m+n,$ $r \in \Z/p\Z,$ and $\lambda \in X(T),$ $i$ is $r$-good for $\lambda$ if and only if $i$ is $r$-normal for $\lambda$ and $i$ is $r$-conormal for $\lambda-\varepsilon_{i}.$
\end{lemma}

The following lemma links the notions of normal and good with conormal and cogood, respectively.  Let us fix some notation which we will use again in subsection~\ref{S:returntonormalandgood}.  Let $w_{0}$ denote the longest element of $S_{m+n}.$  Let $\widetilde{V}$ denote a superspace with $\dim_{k}\widetilde{V}_{\0 }=n$ and $\dim_{k}\widetilde{V}_{\1 }=m$ and let $\widetilde{G}$ denote $GL(\widetilde{V}).$  Fix a homogeneous basis $\widetilde{v}_{1}, \dotsc ., \widetilde{v}_{m+n}$ satisfying 
\[
\p{\widetilde{v}}_{i}=\p{v}_{w_{0}i}+\1 
\] for all $1 \leq i \leq m+n,$ where $v_{1}, \dotsc , v_{m+n} $ is our usual fixed basis for $V.$   Let us write $\widetilde{T} \cong T$ for the subgroup of all diagonal matrices of $\widetilde{G}.$  Following the procedure discussed in subsection~\ref{SS:crystalofXofT}, we put a crystal structure on $X(\widetilde{T})=\bigoplus_{i=1}^{m+n}\Z\varepsilon_{i}$ lifted from the dual crystal 
\[
\left( \widetilde{B}=\mathcal{B}_{\p{\widetilde{v}}_{1}} \otimes \dotsb \otimes \mathcal{B}_{\p{\widetilde{v}}_{m+n}}, \tilde{e}_{r}^{*}, \tilde{f}_{r}^{*}, \varepsilon_{r}^{*}, \varepsilon_{r}^{*}, \wt \right) .
\]  

Let $\tilde{s}: X(T) \to  X(\widetilde{T})$ be given by 
\begin{equation}\label{E:tildesdef}
\varepsilon_{i} \mapsto -\varepsilon_{w_{0}i}
\end{equation}
for $1 \leq i \leq m+n.$   

\begin{lemma}\label{L:normaltoconormal}  Let $\lambda \in X(T)$ and $1 \leq t \leq m+n.$  Then $t$ is normal for $\lambda$ if and only if $w_{0}t$ is conormal for $\tilde{s}(\lambda).$  Moreover, $t$ is good for $\lambda$ if and only if $w_{0}t$ is cogood for $\tilde{s}(\lambda).$
\end{lemma}

\begin{proof} Let $r \in \Z/p\Z.$  Recall our notation $\sigma_{r}(\lambda)=(\sigma_{r}(\lambda)_{1}, \dots , \sigma_{r}(\lambda)_{m+n})$ for the $r$-signature of $\lambda.$  Given $1 \leq i< j \leq m+n, $ we call $(i,j)$ an \emph{$r$-canceling pair} for $\lambda$ if $\sigma_{r}(\lambda)_{i}=-,$ $\sigma_{r}(\lambda)_{j}=+,$ and $\sigma_{k}(\lambda)_{t}=0$ for $i < k< j.$

Given $\lambda \in X(T)$ and $1 \leq t \leq m+n,$ let $\xi$ denote $\sum_{k =1}^{m+n} (-1)^{\p{v}_{k}},$  $r=r_{t}(\lambda),$ and $r'=r - \xi.$ It is straightforward to verify that $(i,j)$ is an $r$-canceling pair for $\lambda$ if and only if $(w_{0}j,w_{0}i)$ is an $r'$-canceling pair for $\tilde{s}(\lambda).$ This along with our combinatorial description of the crystal structure on $X(T)$ immediately implies $t$ is normal for $\lambda$ if and only if $w_{0}t$ is conormal for $\tilde{s}(\lambda).$  It is clear from the definitions that $i$ is good for $\lambda$ if and only if $w_{0}i$ is cogood for $\tilde{s}(\lambda).$
\end{proof}

\subsection{Lowering Operators}\label{S:Def} 
In what follows for $1 \leq i <j \leq m+n$ and $1 \leq k \leq m+n$ we write $E_{i,j}=e_{i,j},$ $F_{i,j}=e_{j,i},$ and $H_{k}=e_{k,k}.$  For $1 \leq j  \leq m+n$, let 
\begin{align}\label{E:murphyoperators}
L_{j} &= \sum_{i=1}^{j-1} (-1)^{\p{v}_{i}}F_{i,j}E_{i,j}. 
\end{align} 
A straightforward calculation shows that the $L_{j}$'s commute with one another and, since they are of weight $0,$ with any element of $\Dist{T}$.

For $1 \leq i, j \leq m+n$ and $1 \leq k \leq m+n-1,$ let
\begin{align}\label{E:cdef}
c_{i,j} = (-1)^{\p{v}_{i}}H_{i}-(-1)^{\p{v}_{j}}H_{j}+(-1)^{\p{v}_{i}}\vartheta_{i}-(-1)^{\p{v}_{j}}\vartheta_{j},
\end{align}
\begin{align}\label{E:ctildedef}
\tilde{c}_{i,j} = c_{i,j} + (-1)^{\p{v}_{i}},
\end{align} and let
\begin{align}\label{E:bdef} b_{i,k} = (-1)^{\p{v}_{i}}H_{i}-(-1)^{\p{v}_{k+1}}H_{k+1}+(-1)^{\p{v}_{i}}\vartheta_{i}-(-1)^{\p{v}_{k}}\vartheta_{k}.
\end{align} Note that this notation is compatible with that of subsection~\ref{S:combinatorialdefinitions} in the sense that if $M$ is a $\Dist{G}$-supermodule with $v\in M_{\lambda}$  then the elements $c_{i,j}, b_{i,k} \in \Dist{G}$ act on $v$ by the scalars $c_{i,j}(\lambda)$ and $b_{i,k}(\lambda)$, respectively. 

For $1 \leq i < j \leq m+n$  and $A \subseteq (i..j)$, define
\begin{align}\label{E:loweringopdef}
\widetilde{S}_{i,j}(A) &=\left(\prod_{t \in A} (\tilde{c}_{i,t} - L_{i+1} - \dotsb - L_{t}) \right)F_{i,j}.
\end{align}
Note that our observations on the commutativity of the $L_{j}$'s imply that the order of the product does not need to be specified.  

Let $B^{-}$ denote the opposite Borel subgroup of $G$ consisting of all lower triangular invertible matrices of the form \eqref{E:matrix}.  Let $U$ denote the unipotent radical of $B$ given by letting $U(A) \subseteq B(A)$ be the set of upper triangular unipotent matrices for any commutative superalgebra $A$.  Fix an ordering of the PBW basis (see Lemma~\ref{stbas}) of $\Dist{G}$ so that each PBW monomial is of the form $XY$ where $X \in \Dist{B^{-}}$ and $Y \in \Dist{U}.$   We define the \emph{lowering operators} $S_{i,j}(A)$ by expanding $\widetilde{S}_{i,j}(A)$ in the PBW basis given by this ordering and taking $S_{i,j}(A)$  to be the sum of those terms lying in $\Dist{B^{-}}.$   The idea is that we are interested applying the lowering operators to primitive vectors and the nonconstant elements of $\Dist{U}$ annihilate all primitive vectors.   Consequently, in many of our calculations we will work modulo a left ideal of $\Dist{G}$ which annihilates any primitive vector.   The following lemma illustrates this point of view.

\pagebreak
\begin{lemma}\label{L:techlemma1}  Let $\mathfrak{J}$ be the left ideal of $\Dist{G}$ generated by the nonconstant elements of $\Dist{U}.$  Let $1 \leq i <t<j \leq m+n.$   Then:
\begin{enumerate}
\item [(i)] $L_{t}F_{i,j} \equiv -aF_{i,t}F_{t,j} \: (mod \: \mathfrak{J})$ where $a=(-1)^{\p{v}_{i}+(\p{v}_{i}+\p{v}_{t})(\p{v}_{i}+\p{v}_{j})}$;
\item [(ii)] $(\tilde{c}_{i,t}-L_{i+1}-\dotsb -L_{t})L_{t}F_{i,j} \equiv 0 \: (mod \: \mathfrak{J})$;
\end{enumerate}
\end{lemma}

\begin{proof} Throughout we write $\equiv$ for congruent modulo $\mathfrak{J}$.  

To prove $(i)$ we simply calculate the left hand side working modulo $\mathfrak{J}$:
\begin{align*}
L_{t}F_{i,j} &= \sum_{r < t}(-1)^{\p{v}_{r}}e_{t,r}e_{r,t}e_{j,i}\\
             &= \sum_{r < t}(-1)^{\p{v}_{r}}e_{t,r}((-1)^{(\p{v}_{i}+\p{v}_{j})(\p{v}_{r}+\p{v}_{t})}e_{j,i}e_{r,t}+[e_{r,t}, e_{j,i}])\\
             &\equiv -\sum_{r < t}(-1)^{\p{v}_{r}+(\p{v}_{r}+\p{v}_{t})(\p{v}_{j}+\p{v}_{i})}\delta_{i,r}e_{t,r}e_{j,t}\\
             &= -(-1)^{\p{v}_{i} +(\p{v}_{i}+\p{v}_{t})(\p{v}_{i}+\p{v}_{j})}e_{t,i}e_{j,t}\\
             &= -aF_{i,t}F_{t,j}.
\end{align*}

To prove $(ii)$ one uses $(i)$ and calculates that 
\begin{align*}
(\tilde{c}_{i,t}-L_{i+1}-\dotsb -L_{t-1})L_{t}F_{i,j} &\equiv -a \tilde{c}_{i,t}F_{i,t}F_{t,j} + a\sum_{i < r < t} L_{r}F_{i,t}F_{t,j} \\
&\equiv -a \tilde{c}_{i,t}F_{i,t}F_{t,j} \\
& \hspace{.5in}- \sum_{i < r < t} (-1)^{(\p{v}_{j}+\p{v}_{r})(\p{v}_{i}+\p{v}_{t})}F_{i,r}F_{r,t}F_{t,j}.
\end{align*}
On the other hand, one calculates that 
\begin{align*}
L_{t}L_{t}F_{i,j} &\equiv -a L_{t}F_{i,t}F_{t,j} \\
&\equiv -a \sum_{r < t} (-1)^{\p{v}_{r}}e_{t,r}e_{r,t}e_{t,i}e_{j,t}.
\end{align*}
However, 
\[
e_{t,r}e_{r,t}e_{t,i}e_{j,t} \equiv \begin{cases} 0 &\mbox{ if } r <i; \\
\left( e_{i,i}-(-1)^{\p{v}_{i}+\p{v}_{t}}e_{t,t}+(-1)^{\p{v}_{i}+\p{v}_{t}}+1 \right)e_{t,i}e_{j,t} &\mbox{ if } r=i; \\
e_{t,i}e_{j,t}+(-1)^{(\p{v}_{t}+\p{v}_{r})(\p{v}_{r}+\p{v}_{i})}e_{r,i}e_{t,r}e_{j,t} &\mbox{ if } r > i.
\end{cases}
\]  Together, these imply 
\[
L_{t}L_{t}F_{i,j} \equiv -a \tilde{c}_{i,t}F_{i,t}F_{t,j}-a \sum_{i < r < t} (-1)^{(\p{v}_{t}+\p{v}_{r})(\p{v}_{r}+\p{v}_{i})+\p{v}_{r}}F_{i,r}F_{r,t}F_{t,j}.
\]  Finally, through the miracle of $\mathbb{Z}_{2}$, we have $\p{v}_{i}+(\p{v}_{i}+\p{v}_{t})(\p{v}_{i}+\p{v}_{j})+(\p{v}_{t}+\p{v}_{r})(\p{v}_{r}+\p{v}_{i})+\p{v}_{r}=(\p{v}_{j}+\p{v}_{r})(\p{v}_{i}+\p{v}_{t})$ which implies the desired result.
\end{proof}

A useful fact about the lowering operators, both for calculating them and for proving results, is the recurrence relation given in the following theorem.

\begin{theorem}\label{T:recurrancerelation}  Let $i \leq i < j \leq m+n,$ and let $A \subseteq (i..j)$.  If $A = \emptyset$, then $S_{i,j}(A)=F_{i,j}$.  Otherwise, take any $k \in A$ and let $h= \max([i..k)\backslash A),$ then we have 
\begin{align}\label{E:recurrancerelation}
S_{i,j}(A) &= S_{i,j}(A\backslash \{k \})c_{h,k} + (1-\delta_{h,i})S_{i,j}(A\backslash \{k \} \cup \{h \}) + aS_{i,k}(A_{i..k})S_{k,j}(A_{k..j}),  
\end{align} where $a=(-1)^{\p{v}_{i}+(\p{v}_{i}+\p{v}_{k})(\p{v}_{i}+\p{v}_{j})}$.
\end{theorem}

\begin{proof} Throughout we write $\equiv$ for equivalence modulo the left ideal $\mathfrak{J}$ generated by the nonconstant elements of $\Dist{U}.$  
 If $h \neq i$, then 
\begin{equation*}
\tilde{c}_{i,k}-L_{i+1}-\dotsb -L_{k} = c_{h,k} + \left(\tilde{c}_{i,h}-L_{i+1}-\dotsb -L_{h}\right) -\left(L_{h+1}+\dotsb +L_{k-1} \right) - L_{k}.
\end{equation*}
Using this equation, we replace the $k$th term in the product in the definition of $\widetilde{S}_{i,j}(A)$.   Distributing yields $S_{i,j}(A) \equiv \widetilde{S}_{i,j}(A) \equiv P_{1}+P_{2}+P_{3}+P_{4}$ where
\begin{align*} 
P_{1}&=\prod_{t \in A_{i..k}} \left(\tilde{c}_{i,t}-L_{i+1}-\dotsb -L_{t} \right)c_{h,k}\prod_{t \in A_{k..j}} \left(\tilde{c}_{i,t}-L_{i+1}-\dotsb -L_{t} \right)F_{i,j} \\
&\equiv S_{i,j}(A\backslash \{k \})c_{h,k}
\intertext{ since, as $i < h,k < j$, the $c_{h,k}$ commutes with all the terms;}
P_{2}&=\prod_{t \in A_{i..k}} \left(\tilde{c}_{i,t}-L_{i+1}-\dotsb -L_{t} \right) \left( \tilde{c}_{i,h}-L_{i+1}-\dotsb -L_{h} \right)\\
& \hspace{2.5in} \cdot \prod_{t \in A_{k..j}} \left(\tilde{c}_{i,t}-L_{i+1}-\dotsb -L_{t} \right) F_{i,j}\\
&\equiv S_{i,j}(A\backslash \{k \} \cup \{h \}); \\
P_{3}&=-\prod_{t \in A_{i..k}} \left(\tilde{c}_{i,t}-L_{i+1}-\dotsb -L_{t} \right) \left( L_{h+1}+\dotsb + L_{k-1}\right) \\
& \hspace{2.5in} \cdot \prod_{t \in A_{k..j}} \left(\tilde{c}_{i,t}-L_{i+1}-\dotsb -L_{t} \right) F_{i,j} \\
&=- \sum_{r=h+1}^{k-1}  \prod_{t \in A_{i..k}} \left(\tilde{c}_{i,t}-L_{i+1}-\dotsb -L_{t} \right)  \prod_{t \in A_{k..j}} \left(\tilde{c}_{i,t}-L_{i+1}-\dotsb -L_{t} \right)L_{r} F_{i,j}\\
&\equiv 0\\
\intertext{by Lemma \ref{L:techlemma1}(ii); and,}
P_{4}&=-\prod_{t \in A_{i..k}} \left(\tilde{c}_{i,t}-L_{i+1}-\dotsb -L_{t} \right)L_{k}\prod_{t \in A_{k..j}} \left(\tilde{c}_{i,t}-L_{i+1}-\dotsb -L_{t} \right)F_{i,j} \\
&\equiv a S_{i,k}(A_{i..k})S_{k,j}(A_{k..j}) 
\end{align*}  by Lemma~\ref{L:techlemma1}(i) and a careful calculation.  Together these imply the result.

The case when $h=i$ is argued similarly using the equation 
\begin{equation*}
\tilde{c}_{i,k}-L_{i+1}-\dotsb -L_{k} = \tilde{c}_{i,k}-\left(L_{i+1}+\dotsb + L_{k-1} \right)-L_{k}.
\end{equation*}
\end{proof}

\subsection{Technical Lemmas}\label{S:technicallemmas}

We now prove several technical lemmas about the lowering operators which will be useful in what follows.  The proof of the lemmas in this subsection follow the arguments of the analogous results of Brundan in \cite{B2}.  The exception is the proof of Lemma~\ref{L:raisingloweredvectors}.  This proof is new and is simpler than Brundan's proof of the analogous result, which requires the introduction of certain formal polynomials.

The following lemma records how the lowering operators commute with $E_{l}=E_{l,l+1}.$  

\pagebreak
\begin{lemma}\label{L:commutatorlemma}  Let $1 \leq i< j \leq m+n$ and let $A \subseteq (i..j)$.  Let $1 \leq l \leq m+n-1$ and let $J_{l}$ be the left ideal of $\dist{}$ generated by $E_{l}$.  Then,
\begin{enumerate}
\item [(i)]  If either \begin{enumerate}
             \item [(a)] $l+1 \in A$, or
             \item [(b)] $l \notin \{i \} \cup A$ and $l+1 \notin A \cup \{j \}$.
             \end{enumerate} then 
\begin{equation}\label{}
E_{l}S_{i,j}(A) \equiv 0 \: (mod \: J_{l});
\end{equation}
\item [(ii)] If $l \in \{i \} \cup A$ and $l+1 \notin A \cup \{j \}$, then 
\[
E_{l}S_{i,j}(A) \equiv -b S_{i,l}(A_{i..l})S_{l+1,j}(A_{l+1..j})  \: (mod \: J_{l}),
\] where 
\[
b= \begin{cases} (-1)^{(\p{v}_{i}+\p{v}_{i+1})(\p{v}_{i}+\p{v}_{j})}, & \mbox{ if } l=i; \\
(-1)^{\p{v}_{i}+(\p{v}_{i}+\p{v}_{l+1})(\p{v}_{i}+\p{v}_{j})}, & \mbox{ if } l \neq i.
\end{cases}
\]
\item [(iii)] If $l = j-1 \notin \{i \}\cup A$, then 
\[
E_{l}S_{i,j}(A) \equiv S_{i,j-1}(A)   \: (mod \: J_{l}).
\] 
\item [(iv)] If $l = j-1 \in A$ and $h = \max([i\dotsb l)\backslash A)$, then 
\[
E_{l}S_{i,j}(A) \equiv S_{i,j-1}(A_{i..j-1})b_{h,j-1} + (1-\delta_{i,h})S_{i,j-1}(A \backslash \{j-1 \} \cup \{h \})  \: (mod \: J_{l}).
\] where $b_{h,j-1}$ is as defined in \eqref{E:bdef}.
\end{enumerate}
\end{lemma}

\begin{proof}  Throughout we write $\equiv$ for equivalence modulo $J_{l}.$
 Define the \emph{height}  of a subset $A \subseteq \{1, \dotsc , m+n \}$ by 
\[
\height(A)= \sum_{ k \in A} k.
\]  

To prove $(i)$ we induct on $\height(A)$.  The base cases are if $A= \emptyset$ or if $A \neq \emptyset$ and $\height(A)=\height(\{i+1 \})$, the minimal nonzero value the height function obtains.  If $A=\emptyset$ then necessarily condition $(b)$ holds.   A direct calculation shows that $E_{l}S_{i,j}(A)=E_{l}F_{i,j} \equiv 0$.  If $\height(A)=\height(\{i+1 \})$ then $A=\{i+1 \}$ and by Theorem \ref{T:recurrancerelation} we have 
\begin{equation*}
S_{i,j}(A)= F_{i,j}c_{i,i+1} + (-1)^{\p{v}_{i}+(\p{v}_{i}+\p{v}_{i+1})(\p{v}_{i}+\p{v}_{j})}F_{i,i+1}F_{i+1,j}.
\end{equation*}  A calculation using this and condition $(b)$ shows $E_{l}S_{i,j}(A) \equiv 0$.  If $A=\{i+1 \}$ and condition $(a)$ holds, then $l=i$ and again a calculation shows $E_{l}S_{i,j}(A) \equiv 0$.

Now suppose $\height(A) > \height(\{i+1 \})$ and the result holds for all sets of smaller height.  If $A=\{l+1 \}$, then by our inductive assumption $l \neq i$ and by applying Theorem \ref{T:recurrancerelation} twice we obtain 
\begin{multline*}
S_{i,j}(A)=F_{i,j}c_{l-1,l+1}+(1-\delta_{i,l-1})S_{i,j}(\{l-1 \}) + (-1)^{\p{v}_{i}+(\p{v}_{i}+\p{v}_{l})(\p{v}_{i}+\p{v}_{j})}F_{i,l}F_{l,j}\\
             +(-1)^{\p{v}_{i}+(\p{v}_{i}+\p{l+1})(\p{v}_{i}+\p{v}_{j})}F_{i,l+1}F_{l+1,j}.
\end{multline*}
The result then follows from a direct calculation and the observation that $\height(\{l-1 \}) = \height(A)-2$ so by the inductive hypothesis $E_{l}S_{i,j}(\{l-1 \}) \equiv 0$.  Otherwise, if $A \neq \{l+1 \}$, then we can choose $k \in A$ with $k \neq l+1$ and apply theorem \ref{T:recurrancerelation} and the inductive hypothesis to prove the result if either condition $(a)$ or $(b)$ holds.

Now we prove $(ii).$  First one proves the case $l=i$ by inducting on $\height(A)$.  If $\height(A)=0$ then $A=\emptyset$ and the result follows from a direct computation.  Now if we assume $\height(A) > 0$ and the statement holds for all lesser heights, then we can choose $k \in  A$ and the result follows from Theorem \ref{T:recurrancerelation}, Lemma \ref{L:commutatorlemma}(i), and the inductive hypothesis.   To prove the case when $l \neq i$, one uses Theorem \ref{T:recurrancerelation} with $k=l$ and the previous case.

One proves $(iii)$ by using Theorem \ref{T:recurrancerelation} to induct on $\height(A)$.

One proves $(iv)$ by applying Theorem \ref{T:recurrancerelation} and the previous parts of the lemma.
\end{proof}

\begin{lemma}\label{L:raisingloweredvectors}  Let $M$ be a $\Dist{G}$-supermodule and let $v \in M_{\lambda}$ be a primitive vector of weight $\lambda \in X(T).$
Let $1 \leq i <j \leq m+n$ and let  $A,B \subseteq (i..j)$ with $|A|=|B|$ and $A \downarrow B$.  Furthermore, assume $c_{i,h}(\lambda) \equiv 0 \modp$ for $h \in (i..j) \backslash A$, and $b_{i,h}(\lambda)\equiv 0 \modp$ for $h \in (i..j) \backslash B$.  Then, 
\begin{equation}\label{E:bequation}
E_{i}\dotsb E_{j-1}S_{i,j}(A).v_{\lambda} = \varepsilon \prod_{t \in \{i \} \cup B}b_{i,t}(\lambda)v_{\lambda}
\end{equation} where $\varepsilon = \pm 1$.
\end{lemma}

\begin{proof}   We prove the statement by induction on $d(i,j) := | (i..j) |$.  The base case is $d(i,j) = 0$.  Then $A=B=\emptyset,$ $j=i+1,$ and the left hand side of \eqref{E:bequation} is 
\begin{align*}
E_{i}F_{i,j}.v_{\lambda} &= e_{i,j}e_{j,i}.v_{\lambda}\\
&= \left(e_{i,i}-(-1)^{\p{v}_{i}+\p{v}_{j}}e_{j,j}\right).v_{\lambda}\\
&= (-1)^{\p{v}_{i}}b_{i,i}(\lambda)v_{\lambda},
\end{align*} which is the right hand side of \eqref{E:bequation}.  In this case, incidentally, $\varepsilon = (-1)^{\p{v}_{i}}$.

Now we assume $d(i,j) > 0$ and that \eqref{E:bequation} holds for all smaller $d(i,j)$.  We proceed by considering cases:

\textbf{Case 1:} Say $j-1 \notin A$.  Observe that this along with $A \downarrow B$ implies $j-1 \notin B$.  Then the left hand side of \eqref{E:bequation} is

\begin{align*}
E_{i}\dotsb E_{j-1}S_{i,j}(A).v_{\lambda} &= E_{i}\dotsb E_{j-2}S_{i,j-1}(A).v_{\lambda} \\
&= \varepsilon \prod_{t \in \{i \} \cup B}b_{i,t}(\lambda)v_{\lambda}
\end{align*} by Lemma~\ref{L:commutatorlemma}(iii) and the inductive assumption. 

\textbf{Case 2:} Say $j-1 \in A$ and $j-1 \in B$.  Set $A_{1}:=A\backslash \{j-1 \}$ and $B_{1}:=B \backslash \{j-1 \}$.  Then $A_{1}, B_{1} \subseteq (i..j-1)$, $|A_{1}| = |B_{1}|$, and $A_{1} \downarrow B_{1}$.  If $h=\max([i..j-1)\backslash A) = i$, then by Lemma~\ref{L:commutatorlemma}(iv) and the inductive assumption the left hand side of \eqref{E:bequation} is 
\begin{align*}
E_{i}\dotsb E_{j-1}S_{i,j}(A).v_{\lambda} &= E_{i}\dotsb E_{j-2}S_{i,j-1}(A_{1})b_{i,j-1}.v_{\lambda}\\
&= b_{i,j-1}(\lambda)E_{i}\dotsb E_{j-2}S_{i,j-1}(A_{1}).v_{\lambda}\\
&=  b_{i,j-1}(\lambda)\varepsilon\left( \prod_{t \in \{i \} \cup B_{1}}b_{i,t}(\lambda) \right)v_{\lambda}\\
&= \varepsilon \prod_{t \in \{i \} \cup B} b_{i,t}(\lambda)v_{\lambda}.
\end{align*} 
On the other hand, if $h > i$, set $A_{2}:=A_{1} \cup \{h \}$ and $B_{2}:=B_{1} \cup \{l \}$, where $l :=\min((i..j) \backslash B)$ (Note that $l$ exists as $h > i$ so $A \neq (i..j)$ hence, since $|A|=|B|$, $B \neq (i..j)$).  Then $A_{2}, B_{2} \subseteq (i..j-1)$, $|A_{2}|=|B_{2}|$, and $A_{2} \downarrow B_{2}$.  Again, by Lemma~\ref{L:commutatorlemma}(iv) and the inductive assumption the left hand side of \eqref{E:bequation} is 
\begin{align*}
E_{i}\dotsb E_{j-1}S_{i,j}(A).v_{\lambda} &= E_{i}\dotsb E_{j-2}S_{i,j-1}(A_{1})b_{h,j-1}.v_{\lambda} + E_{i}\dotsb E_{j-2}S_{i,j-1}(A_{2}).v_{\lambda}\\
&=b_{h,j-1}(\lambda)\varepsilon \left(\prod_{t \in \{i \} \cup B_{1}}b_{i,t}(\lambda)\right)v_{\lambda}+\varepsilon' \prod_{t \in \{i \} \cup B_{2}} b_{i,t}(\lambda)v_{\lambda} 
\intertext{However, $h \in (i..j) \backslash A$ so $c_{i,h}(\lambda)= 0$, hence $b_{h,j-1}(\lambda)=b_{i,j-1}(\lambda)-c_{i,h}(\lambda) = b_{i,j-1}(\lambda)$.  Also, $l \in (i..j) \backslash B$ so $b_{i,l}(\lambda)= 0$.  Together, these yield}
&= \varepsilon \prod_{t \in \{i \} \cup B}b_{i,t}(\lambda)v_{\lambda}
\end{align*}
\textbf{Case 3:} Say $j-1 \in A$ and $j-1 \notin B$.  Since $j-1 \notin B$ and $|A|=|B|$, we see that $h=\max([i..j-1) \backslash A) \neq i$.  Note, too, that $A\backslash \{j-1 \}\cup \{h \} \downarrow B$.  Consequently, by Lemma~\ref{L:commutatorlemma}(iv) and the inductive assumption, the left hand side of \eqref{E:bequation} is 
\begin{align*}
E_{i}\dotsb E_{j-1}S_{i,j}(A).v_{\lambda} &= E_{i}\dotsb E_{j-2}S_{i,j-1}(A \backslash \{j-1 \})b_{h,j-1}.v_{\lambda} \\
&\hspace{1in} +  E_{i}\dotsb E_{j-2}S_{i,j-1}(A\backslash \{j-1 \} \cup \{h \}).v_{\lambda} \\
&= b_{h,j-1}(\lambda)E_{i}\dotsb E_{j-2}S_{i,j-1}(A \backslash \{j-1 \}).v_{\lambda} + \varepsilon \prod_{t \in \{i \}\cup B}b_{i,t}(\lambda)v_{\lambda}
\intertext{However, since $h \in (i..j) \backslash A$ and $j-1 \in (i..j) \backslash B$, we have $b_{h,j-1}(\lambda)=b_{i,j-1}(\lambda)-c_{i,h}(\lambda) = 0$.  Therefore,}
&= \varepsilon \prod_{t \in \{i \} \cup B}b_{i,t}(\lambda).
\end{align*}
Therefore, under the assumptions given \eqref{E:bequation} always holds. 
\end{proof}

\begin{lemma}\label{L:combinatorialcriteria}  Let $1 \leq i < j \leq m+n$, $\lambda \in X(T),$ and $v_{\lambda} \in M(\lambda)_{\lambda}.$  Let $C \subseteq C_{i,j}(\lambda)$ and set $A=(i..j) \backslash C$.  If $S_{i,j}(A).v_{\lambda} \notin \Rad M(\lambda)$, then $B_{i,j}(\lambda) \downarrow C$. 
\end{lemma}

\begin{proof}  We again induct on $d(i,j) := |(i..j)|$.  If $d(i,j)=0$, then $A=C=\emptyset$ and $j=i+1$.  Consequently, we have 
\[
S_{i,j}(A).v_{\lambda}=F_{i,i+1}.v_{\lambda} = e_{i+1,i}.v_{\lambda} \notin \Rad M(\lambda) .
\] However, this implies we have 
\begin{align*} e_{i,i+1}e_{i+1,i}.v_{\lambda}&=\left(e_{i,i}-(-1)^{\p{v}_{i}+\p{v}_{i+1}}e_{i+1,i+1}+(-1)^{\p{v}_{i}+\p{v}_{i+1}}e_{i+1,i}e_{i,i+1} \right).v_{\lambda} \\
&=\left(e_{i,i}-(-1)^{\p{v}_{i}+\p{v}_{i+1}}e_{i+1,i+1} \right).v_{\lambda} \\
&=(-1)^{\p{v}_{i}}b_{i,i}.v_{\lambda} \\
& \notin \Rad M(\lambda).
\end{align*}  In particular, this implies $b_{i,i}.v_{\lambda} \neq 0$, hence $b_{i,i}(\lambda) \neq 0$.  Therefore $B_{i,j}(\lambda) = \emptyset$ and, trivially, $B_{i,j}(\lambda) \downarrow C$.

Now assume $d(i,j) > 0$ and that the statement holds for all smaller $d(i,j)$.  First, we observe that $E_{l}S_{i,j}(A).v_{\lambda} \notin \Rad  M(\lambda)$ for some $1 \leq l \leq m+n-1 $.  Otherwise this along with a weight argument implies $S_{i,j}(A).v_{\lambda} \in \Rad  M(\lambda)$, a contradiction.  Thus we can choose an $1 \leq l \leq m+n-1 $ with $E_{l}S_{i,j}(A).v_{\lambda} \notin \Rad M(\lambda)$.  We now consider the possibilities for $l:$

\textbf{Case 1:} $l \in \{i \}\cup A, l+1 \notin A \cup \{j \}$.  Then by Lemma \ref{L:commutatorlemma}(ii) we have 
\[
E_{l}S_{i,j}(A).v_{\lambda}= -bS_{i,l}(A_{i..l})S_{l+1,j}(A_{l+1..j}).v_{\lambda} \notin \Rad M(\lambda)
\] where $b$ is as in the Lemma.  However $S_{i,l}(A_{i..l})$ and $S_{l+1,j}(A_{l+1..j})$ commute so both $S_{i,l}(A_{i..l}).v_{\lambda} \notin \Rad  M(\lambda)$ and $S_{l+1,j}(A_{l+1..j}).v_{\lambda} \notin \Rad  M(\lambda)$.  By induction $B_{i,l}(\lambda) \downarrow C_{i..l} := (i..l) \backslash A_{i..l}$ and $B_{l+1, j} \downarrow C_{l+1..j}:= (l+1..j) \backslash A_{l+1..j}$.  However, in the case under consideration, $l+1 \in C_{i,j}(\lambda)$.  That is, $c_{i, l+1}(\lambda) = 0$ so $b_{i, l}(\lambda) \neq 0$, hence $l \notin B_{i,j}(\lambda)$.  We also observe that $B_{l+1,j}(\lambda) \subseteq B_{i,j}(\lambda)$ and $B_{i,l}(\lambda) \subseteq B_{i,j}(\lambda)$.  Thus we have $B_{i,l}(\lambda) \cup B_{l+1, j}(\lambda) = B_{i,j}(\lambda)$ (with equality coming from the fact that $l \notin B_{i,j}(\lambda)$).  However, $C=C_{i..l}\cup C_{l+1..j} \cup \{l+1 \}$.  Therefore $B_{i,j}(\lambda) \downarrow C$.

\textbf{Case 2:} $l = j-1$, $l \notin \{i \} \cup A$.  By Lemma \ref{L:commutatorlemma}(iii) we have 
\[
E_{l}S_{i,j}(A).v_{\lambda} = S_{i,j-1}(A).v_{\lambda} \notin \Rad  M(\lambda).
\]  By induction, we have $B_{i,j-1}(\lambda) \downarrow C_{i..j-1}:=(i..j-1) \backslash A$.  But $C = C_{i..j-1} \cup \{j-1 \}$ and either $B_{i,j}(\lambda)= B_{i,j-1}(\lambda)$ or $B_{i,j}(\lambda)=B_{i,j-1}(\lambda) \cup \{j-1 \}.$  In either case, the weakly increasing injective map $\theta : B_{i,j-1}(\lambda) \hookrightarrow C_{i..j-1}$ extends to a weakly increasing injective map $B_{i,j}(\lambda) \hookrightarrow C$.  Therefore $B_{i,j}(\lambda) \downarrow C$.

\textbf{Case 3:} $l = j-1 \in A$.  Let $h= \max([i..l) \backslash A)$.  Then by Lemma \ref{L:commutatorlemma} (iv)
\[
E_{l}S_{i,j}(A).v_{\lambda} = S_{i,j-1}(A_{i..j-1})b_{h,j-1}.v_{\lambda} + (1-\delta_{i,h})S_{i,j-1}(A \backslash \{j-1 \} \cup \{h \}).v_{\lambda}.
\]  Then necessarily either $S_{i,j-1}(A_{i..j-1})b_{h,j-1}.v_{\lambda} \notin \Rad  M(\lambda)$ or $(1-\delta_{i,h})S_{i,j-1}(A \backslash \{j-1 \} \cup \{h \}).v_{\lambda} \notin \Rad  M(\lambda)$.  If the former occurs, then by induction we have $B_{i,j-1}(\lambda) \downarrow C_{i..j-1}$.  As in Case 2, one concludes that $B_{i,j}(\lambda) \downarrow C$.  If the latter occurs, then $h \neq i$ and by induction $B_{i,j-1}(\lambda) \downarrow C \backslash \{h \}$.  Again one can extend any weakly increasing injective map to prove that $B_{i,j}(\lambda) \downarrow C$.

All other possibilities are eliminated by Lemma~\ref{L:commutatorlemma} and our assumption that $E_{l}S_{i,j}(A).v_{\lambda} \not\in \Rad M(\lambda).$  This proves the desired result in all possible situations.
\end{proof}

\subsection{Filtrations and Hom-spaces}\label{SS:filtrations}

By weights and Frobenious reciprocity, we have the following lemma.

\begin{lemma}\label{L:homspaces}  Let $\lambda, \mu \in X(T)$ and let $W$ be a nonzero submodule of $W(\mu).$  Then 
\[
\dim_{k} \Hom_{\Dist{G}}(M(\lambda), W)=\delta_{\lambda, \mu}.
\]
\end{lemma}

\begin{corollary}\label{C:countinghoms}  Say $M$ is a $\Dist{G}$-supermodule with a filtration by co-Verma supermodules:
\[
0 = M_{0} \subseteq M_{1} \subseteq M_{2} \subseteq \dotsb \subseteq M_{r}=M
\] with $M_{i}/M_{i-1}$ isomorphic to $W(\mu^{(i)})$ for $1 \leq i \leq r$.  Then 
\[
\dim_{k} \Hom_{\Dist{G}}(M(\lambda), M) = |\{i: \mu^{(i)}=\lambda \}|.
\] 

\end{corollary}

\begin{proof}  We prove the theorem by inducting on $r.$  The base case of $r=1$ follows immediately from Lemma~\ref{L:homspaces}.  One proves the inductive step by applying the functor $\Hom_{\Dist{G}}(M(\lambda), -)$ to the short exact sequence 
\[
0 \to W(\mu^{(1)}) \to M \to M/M_{1} \to 0,
\] and Lemmas~\ref{L:exttriviality} and~\ref{L:homspaces}.
\end{proof}

For $1 \leq i \leq m+n,$ let $w_{i} \in V^{*}$ denote the element defined by $w_{i}(v_{j})=(-1)^{\p{v}_{i}}\delta_{i,j}.$  Then $\{w_{i}: 1 \leq i \leq m+n \}$ forms a basis for $V^{*}.$

\begin{lemma}\label{L:filtrationofMotimesVdual}  Let $\lambda \in X(T).$  Then we have the filtration
\[
0 =M_{0}\subseteq M_{1} \subseteq \dotsb \subseteq M_{m+n}=M(\lambda)\otimes V^{*}
\] of $M(\lambda)\otimes V^{*}$ with $M_{i}/M_{i-1} \cong M(\lambda - \varepsilon_{m+n-i+1})$ for all $i=1, \dotsc , m+n.$  Furthermore, if $v_{\lambda}$ is the canonical generator for $M(\lambda),$ then the image of $v_{\lambda} \otimes w_{m+n-i+1}$ is the canonical generator of $M_{i}/M_{i-1}.$
\end{lemma}

\begin{proof} Let $k_{\lambda}$ be the $\Dist{B}$-supermodule of highest weight $\lambda \in X(T),$ and let $v_{\lambda}$ denote its canonical generator.  We have a filtration of $k_{\lambda} \otimes V^{*}$ as a $\Dist{B}$-supermodule given by 
\[
0=Q_{0} \subseteq Q_{1} \subseteq \dotsb \subseteq Q_{m+n}=k_{\lambda} \otimes V^{*},
\] where $Q_{k}$ is generated by $v_{\lambda} \otimes w_{m+n}, \dotsc , v_{\lambda} \otimes w_{m+n-k+1}.$  Let 
\[
M_{k}=\Dist{G} \otimes_{\Dist{B}} Q_{k}.
\]
The result follows by exactness of the functor $\Dist{G} \otimes_{\Dist{B}} -$ and the super version of the Tensor Identity (c.f. \cite[I.3.6]{J1}).
\end{proof}

\begin{corollary}\label{C:filtrationofWotimesVdual} Let $\lambda \in X(T).$  Then 
\begin{equation*}
\Hom_{\Dist{G}}(M(\mu), W(\lambda)\otimes V^{*})= \begin{cases} k, &\text{ if $\mu = \lambda -\varepsilon_{i}$ for some $1 \leq i \leq m+n $;}\\
                                                                0, &\text{ otherwise.}

\end{cases}
\end{equation*}
\end{corollary}

\begin{proof}  It follows from Lemma~\ref{L:filtrationofMotimesVdual} and taking contravariant duals that we have the filtration 
\begin{equation}\label{E:dualfiltration}
0 =N_{0}\subseteq N_{1} \subseteq \dotsb \subseteq N_{m+n}=W(\lambda)\otimes V^{*}
\end{equation}
with $N_{i}/N_{i-1} \cong W(\lambda - \varepsilon_{i})$ for all $i=1, \dotsc , m+n.$ The result then follows by applying Corollary~\ref{C:countinghoms}. 
\end{proof}

\subsection{Return to Normal, Good, Conormal, and Cogood}\label{S:returntonormalandgood}  We are now able to provide representation theoretic interpretations of the combinatorial notions of normal, good, conormal, and cogood.  Note that the arguments used here are an adaptation of those used in \cite{B2}.  

Let $V'$ denote the subspace of $V$ spanned by the vectors $v_{1}, \dotsc , v_{m+n-1}.$  We consider $G'=GL(V')$ as a subgroup of $G=GL(V)$ in the natural way and make the corresponding identification of $\Dist{G'}$ as a subalgebra of $\Dist{G}.$  A direct calculation verifies the following lemma.

\begin{lemma}\label{L:eandfhomomorphism} Let $M$ be a $\Dist{G}$-supermodule.  Define the following linear map 
\[
\begin{array}{ll}\label{E:efdefs}
e: M \to M\otimes V^{*}, & x \mapsto x \otimes w_{m+n} + \sum_{h=1}^{m+n-1} (-1)^{(\p{v}_{h}+\p{x})(\p{v}_{h}+\p{v}_{m+n})}E_{h,m+n}x \otimes w_{h}.
\end{array}
\] The map $e$ is a homomorphism of $\Dist{G'}$-supermodules.
\end{lemma}

\begin{theorem}\label{T:normalitytheorem}  Let $\lambda \in X(T)$ and let $1 \leq i \leq m+n.$  Let $v_{\lambda} \in M(\lambda)_{\lambda}$ denote the canonical generator of $M(\lambda).$  The following are equivalent:
\begin{enumerate}
\item [(i)] $\dim_{k} \Hom_{\Dist{G}}(M(\lambda-\varepsilon_{i}), L(\lambda)\otimes V^{*}) = 1;$
\item [(ii)] $S_{i,m+n}(A).v_{\lambda} \notin \Rad  M(\lambda)$, where $A:= (i..m+n) \backslash C_{i,m+n}(\lambda)$;
\item [(iii)] $B_{i,m+n}(\lambda) \downarrow C_{i,m+n}(\lambda)$ (i.e. $i$ is normal for $\lambda$);
\end{enumerate}
\end{theorem}

\begin{proof} 

Say $i=m+n.$  A consideration of weights shows that the image of $v_{\lambda} \otimes w_{m+n} \in L(\lambda) \otimes V^{*}$ is a primitive vector.  Consequently, by Frobenious reciprocity we have $\Hom_{\Dist{G}}(M(\lambda -\varepsilon_{m+n}), L(\lambda) \otimes V^{*}) \neq 0.$ When $i=m+n$ it is straightforward to verify that $(ii)$ and $(iii)$ also always hold.  Therefore, for the rest of the proof we assume $1 \leq i \leq m+n-1.$

$(i) \Rightarrow (ii):$ Let $\varphi: M(\lambda)\to W(\lambda-\varepsilon_{i})\otimes V$ be the non-zero homomorphism given in Corollary~\ref{C:filtrationofWotimesVdual}.  We can choose $\varphi$ so that $\varphi(v_{\lambda})=w\otimes v_{i},$ where $w$ is the dual of the canonical generator of $M(\lambda-\varepsilon_{i}).$  Observe that $(i)$ holds if and only if $\Rad M(\lambda) \subseteq \Ker \varphi$.  Suppose $(ii)$ is false, then $S_{i,m+n}(A).v_{\lambda}=\widetilde{S}_{i,m+n}(A).v_{\lambda} \in \Rad  M(\lambda)$.  It suffices to show $\widetilde{S}_{i,n}(A).v_{\lambda} \notin \Ker \varphi$.  As a vector space, we have the direct sum decomposition 
\[
W(\lambda-\varepsilon_{i})\otimes V =  \bigoplus_{i=1 }^{m+n}  W(\lambda-\varepsilon_{i})\otimes v_{i}.
\]
From the definition of $\widetilde{S}_{i,m+n}(A),$ one can verify that 
\begin{align*}
\varphi(\widetilde{S}_{i,m+n}(A).v_{\lambda})=\widetilde{S}_{i,m+n}(A).(w\otimes v_{i})&= \prod_{t \in A} \tilde{c}_{i,t}F_{i,m+n}.(w\otimes v_{i}) + (*) \\
&=F_{i,m+n}\prod_{t \in A} c_{i,t}.(w \otimes v_{i}) + (*)\\
&=F_{i,m+n}\prod_{t \in A} c_{i,t}(\lambda)(w \otimes v_{i}) + (*)\\
&=\prod_{t \in A} c_{i,t}(\lambda)(w \otimes v_{m+n} + (F_{i,m+n}w) \otimes v_{i}) + (*)
\end{align*} where $(*)$ lies in $\bigoplus_{i =1}^{m+n-1}  W(\lambda-\varepsilon_{i})\otimes v_{i}.$  In particular, the projection of $\varphi(\widetilde{S}_{i,m+n}(A).v_{\lambda})$ onto $W(\lambda-\varepsilon_{i}) \otimes v_{m+n}$ is precisely $\prod_{t \in A} c_{i,t}(\lambda)(w \otimes v_{m+n}),$ which is nonzero by the choice of $A.$  Therefore $\varphi(\widetilde{S}_{i,m+n}(A).v_{\lambda}) \neq 0.$

$(ii) \Rightarrow (iii):$  This follows from Lemma \ref{L:combinatorialcriteria}.

$(iii) \Rightarrow (i):$  Since $B_{i,m+n}(\lambda) \downarrow C_{i,m+n}(\lambda)$, then we can find a subset $C \subseteq C_{i,m+n}(\lambda)$ such that $B_{i,m+n}(\lambda) \downarrow C$ and $|B_{i,m+n}(\lambda)|=|C|$.   Note that since $B_{i,m+n}(\lambda) \downarrow C_{i,m+n}(\lambda)$, we have $i \notin B_{i,m+n}(\lambda)$.  So, setting $B=(i..m+n) \backslash B_{i,m+n}(\lambda)$ and $A=(i..m+n) \backslash C$, we have $A,B \subseteq (i..m+n)$, $|A|=|B|$, and $A \downarrow B$.  Consequently, by Lemma \ref{L:raisingloweredvectors}, we have 
\[
E_{i}\dotsb E_{m+n-1}S_{i,m+n}(A).v_{\lambda} = \varepsilon \prod_{t \in \{i \} \cup B}b_{i,t}(\lambda)v_{\lambda}.
\] However, by the definition of $B$, the right hand side is a nonzero scalar multiple of $v_{\lambda}$.  Therefore, $z:=S_{i,m+n}(A).v_{\lambda} \notin \Rad  M(\lambda)$ and the image of $z$ in $L(\lambda)$ is nonzero.  Consequently $z \otimes w_{m+n} \in L(\lambda) \otimes V^{*}$ is nonzero and, hence,
\begin{equation}\label{L:defofv(i)}
v_{\lambda}(i):=e(z)=z \otimes w_{m+n} + \sum_{k =1}^{m+n-1} (-1)^{(\p{v}_{k}+\p{v}_{i}+\p{v}_{m+n})(\p{v}_{k}+\p{v}_{m+n})}E_{k,m+n}z \otimes w_{k}
\end{equation} is nonzero.

We now show $E^{(r)}_{l}.z \in \Rad  M(\lambda)$ for all $r \geq 1$ and $ 1 \leq l \leq m+n-2$.  If $r >2$ then $E^{(r)}_{l}.z=0$ by a weight argument.  It remains to consider the case when $r=1:$  If $l +1 \in A,$ or $l \notin \{i \} \cup A$ and $l+1 \notin A$, then $E_{l}.z = 0$ by Lemma \ref{L:commutatorlemma}(i).  If $l \in \{i \} \cup A$ and $l+1 \notin A$, then by Lemma \ref{L:commutatorlemma}(ii), $E_{l}.z=-bS_{i,l}(A_{i..l})S_{l+1..m+n}(A_{l+1..m+n}).v_{\lambda}$.  Since $S_{i,l}(A_{i..l})$ and $S_{l+1,m+n}(A_{l+1..m+n})$ commute, it suffices to show either  $S_{i,l}(A_{i..l}).v_{\lambda} \in \Rad  M(\lambda)$ or $S_{l+1,m+n}(A_{l+1..m+n}).v_{\lambda} \in \Rad  M(\lambda)$.  First, however, observe that our assumptions about $l$ imply $c_{i,l+1}(\lambda) = 0$, so $b_{i,l}(\lambda) \neq 0$, $C_{l+1..m+n}\subseteq C_{l+1,m+n}(\lambda)$, and $B_{i,m+n}(\lambda)=B_{i,l}(\lambda) \cup B_{l+1,m+n}(\lambda)$.  Now suppose $S_{i,l}(A_{i..l}).v_{\lambda} \notin \Rad  M(\lambda)$.  Then by Lemma \ref{L:combinatorialcriteria}, $B_{i,l}(\lambda) \downarrow C_{i..l}$ and, hence, $|B_{i,l}(\lambda)| \leq |C_{i..l}|$.  Thus, $|B_{l+1,m+n}(\lambda)| \geq |C_{l..m+n}|$, so $|B_{l+1, m+n}(\lambda)| > |C_{l+1..m+n}|$.  But this implies that $B_{l+1,m+n}(\lambda) \downarrow C_{l+1..m+n}$ is false.  Therefore, $S_{l+1,m+n}(A_{l+1..m+n}).v_{\lambda} \in \Rad  M(\lambda)$ by Lemma \ref{L:combinatorialcriteria}. 

The preceding paragraph proves $z \in L(\lambda)$ is a $\Dist{G'}$ primitive vector hence, by Lemma~\ref{L:eandfhomomorphism}, $v_{\lambda}(i)$ is a $\Dist{G'}$ primitive vector.  Furthermore, observe that a weight argument along with a direct calculation verifies $E_{m+n-1}^{(r)}.v_{\lambda}(i)=0$ for any $r \geq 1.$  Therefore $v_{\lambda}(i)$ is in fact a $\Dist{G}$ primitive vector.  By a routine application of Frobenious reciprocity we obtain  
\[
\Hom_{\Dist{G}}(M(\lambda-\varepsilon_{i}), L(\lambda) \otimes V^{*}) \neq 0.
\]
\end{proof}

To continue, recall from subsection~\ref{S:combinatorialdefinitions}   the definitions of $\widetilde{V},$  $\widetilde{G}=GL(\widetilde{V}),$ and $X(\widetilde{T}).$   Let us write $\widetilde{B}$ for the Borel subgroup of $\widetilde{G}$ of upper triangular matrices, $\widetilde{L}(\lambda)$ for the irreducible $\Dist{\widetilde{G}}$-supermodule of highest weight $\lambda \in X(\widetilde{T}),$ etc. 

Let $\tilde{S}: \Dist{G} \to \Dist{\widetilde{G}}$ be the even superalgebra isomorphism given by 
\begin{equation}\label{E:Stildedef}
\tilde{S}(e_{i,j}) = -(-1)^{\p{v}_{i}(\p{v}_{i}+\p{v}_{j})}e_{w_{0}j,w_{0}i}
\end{equation}
for all $1 \leq i,j \leq m+n,$ where $w_{0} \in S_{m+n}$ is the longest element.  Given a $\Dist{\widetilde{G}}$-supermodule, $M,$ we can twist it by $\tilde{S}$ to obtain a $\Dist{G}$-supermodule.  Namely, $M=M$ as a superspace with action given by $x.m=\tilde{S}(x)m$ for all $x \in \Dist{G}$ and all $m \in M.$  We denote the resultant $\Dist{G}$-supermodule by $M^{\tilde{S}}.$  Since $\tilde{S}(\Dist{T})=\Dist{\widetilde{T}},$ the weight space decomposition of $M$ is preserved.  In fact, since $e_{i, i} \mapsto -e_{w_{0}i,w_{0}i}$ for all $1 \leq i \leq m+n,$ we have  
\[
M_{\tilde{s}(\lambda)}=(M^{\tilde{S}})_{\lambda}
\] for all $\lambda \in X(T),$ where the map $\tilde{s}:X(T) \to X(\widetilde{T})$ is as defined in \eqref{E:tildesdef}.  Furthermore, since  
$\tilde{S}(\Dist{B})=\Dist{\widetilde{B}},$ twisting by $\tilde{S}$ preserves primitive vectors.  Consequently we have
\[\begin{array}{lcl}
\widetilde{M}(\tilde{s}(\lambda))^{\tilde{S}} \cong M(\lambda) &\text{ and } &\widetilde{L}(\tilde{s}(\lambda))^{\tilde{S}} \cong L(\lambda).
\end{array}
\]
In particular, we have that
\[
(\widetilde{V}^{*})^{\tilde{S}}\cong V.
\]
Note that for any $\Dist{\widetilde{G}}$-supermodules $M$ and $N$ we have 
\[
(M \otimes N)^{\tilde{S}} \cong M^{\tilde{S}}\otimes N^{\tilde{S}}.
\]

\begin{corollary}\label{C:conormalitytheorem}  Let $\lambda \in X(T)$ and let $1 \leq i \leq m+n.$   Then $i$ is conormal for $\lambda$ if and only if   
\[
\dim_{k} \Hom_{\Dist{G}}(M(\lambda+\varepsilon_{i}), L(\lambda)\otimes V) = 1.
\]
\end{corollary}

\begin{proof}
By Lemma~\ref{L:normaltoconormal}, $i$ is conormal for $\lambda \in X(T)$ if and only if $w_{0}i$ is normal for $\tilde{s}(\lambda) \in X(\tilde{T}).$   By the above remarks we have the following isomorphisms of $\Hom$-spaces: 
\begin{align*}
&\Hom_{\Dist{\widetilde{G}}}(\widetilde{M}(\tilde{s}(\lambda)-\varepsilon_{w_{0}i}), \widetilde{L}(\tilde{s}(\lambda))\otimes \widetilde{V}^{*}) \\
&\hspace{.5in}\cong \Hom_{\Dist{G}}(\widetilde{M}(\tilde{s}(\lambda)-\varepsilon_{w_{0}i})^{\tilde{S}}, (\widetilde{L}(\tilde{s}(\lambda))\otimes \widetilde{V}^{*})^{\tilde{S}})\\
&\hspace{.5in}\cong\Hom_{\Dist{G}}(M(\lambda+\varepsilon_{i}), L(\lambda)\otimes V).
\end{align*}  The result then follows immediately from Theorem~\ref{T:normalitytheorem}. 
\end{proof}

\begin{corollary}\label{C:goodcriteria}  Let $1 \leq i \leq m+n$.  Then $i$ is good for $\lambda$ if and only if 
\begin{equation}\label{E:goodhomspace}
\dim_{k} \Hom_{\Dist{G}}(L(\lambda-\varepsilon_{i}), L(\lambda) \otimes V^{*}) =1.
\end{equation}
\end{corollary}

\begin{proof}  By Lemma~\ref{L:normalplusconormalequalsgood}, it suffices to show \eqref{E:goodhomspace} holds if and only if $i$ is normal for $\lambda$ and conormal for $\lambda-\varepsilon_{i}.$  That is, by Theorem~\ref{T:normalitytheorem} and Corollary~\ref{C:conormalitytheorem}, that \eqref{E:goodhomspace} holds if and only if 
\begin{align}\label{E:normalconormalhomspaces}
\Hom_{\Dist{G}}(M(\lambda-\varepsilon_{i}), L(\lambda)\otimes V^{*}) &\neq 0\\
\Hom_{\Dist{G}}(L(\lambda-\varepsilon_{i}), W(\lambda)\otimes V^{*}) &\neq 0. \notag
\end{align}

The $\Hom$-space given in \eqref{E:goodhomspace} is naturally a subspace of the $\Hom$-spaces in \eqref{E:normalconormalhomspaces}.  Hence, if \eqref{E:goodhomspace} holds, then \eqref{E:normalconormalhomspaces} holds.   Conversely, assume \eqref{E:normalconormalhomspaces} holds.  By Corollary~\ref{C:filtrationofWotimesVdual} we have that 
\[
\Hom_{\Dist{G}}(M(\lambda-\varepsilon_{i}), W(\lambda)\otimes V^{*})
\] is exactly one dimensional.  Choose a nonzero homomorphism $f:M(\lambda-\varepsilon_{i}) \to W(\lambda)\otimes V^{*}.$  Since \eqref{E:normalconormalhomspaces}  holds, we have $f(M(\lambda - \varepsilon_{i})) \subseteq L(\lambda) \otimes V^{*}$ and $\operatorname{Ker} f \subseteq \Rad M(\lambda-\varepsilon_{i}).$  Consequently, the map $f$ factors through to give a nonzero $\tilde{f}:L(\lambda-\varepsilon_{i}) \to L(\lambda)\otimes V^{*}.$  This implies the desired result. 
\end{proof}

Arguing as in the proof of Corollary~\ref{C:conormalitytheorem} using the map $\tilde{S}$ one obtains the following corollary of the previous result.

\begin{corollary}\label{C:cogoodcriteria}  Let $1 \leq i \leq m+n$.  Then $i$ is cogood for $\lambda$ if and only if 
\begin{equation}\label{E:cogoodhomspace}
\dim_{k} \Hom_{\Dist{G}}(L(\lambda+\varepsilon_{i}), L(\lambda) \otimes V) =1.
\end{equation}
\end{corollary}

\section{Translation Functors}\label{S:translationfunctors}

\subsection{Translation Functors and Irreducible Supermodules}\label{SS:irreps}  Recall from \eqref{E:transfunctordefs} the definition of the translation functors $E_{r}$ and $F_{r}$ for $r \in \Z/p\Z.$
We can now describe how these functors act on irreducible $\Dist{G}$-supermodules.

\begin{lemma}\label{L:homsintoErLlambda}  Let $\lambda \in X(T)$ and $r \in \Z/p\Z.$ Then there is a filtration
\[0=Q_{0} \subseteq \dots \subseteq Q_{l} = E_{r}L(\lambda)\] 
with $Q_{k}/Q_{k-1}$ a nonzero submodule of $W(\lambda - \varepsilon_{i_{k}})$ and where $i_{1} < \dots < i_{l}$ and $\{i_{1}, \dots , i_{l} \}= \{i : 1 \leq i \leq m+n, i \text{ is $r$-normal for $\lambda$}  \}.$
\end{lemma}

\begin{proof}   First, let $\mu \in X(T).$  We observe that by Theorem~\ref{T:normalitytheorem}
\begin{equation}\label{E:Eahomspace}
\Hom_{\Dist{G}}(M(\mu), E_{r}L(\lambda)) \neq 0
\end{equation}
if and only if $\mu=\lambda -\varepsilon_{i}$ for some $1 \leq i \leq m+n$ which is normal for $\lambda$ and $\wt(\lambda -\varepsilon_{i})=\wt(\lambda)+\gamma_{r}-\gamma_{r+1}.$  However, $\wt(\lambda-\varepsilon_{i}) = \wt(\lambda)-(-1)^{\p{v}_{i}}\gamma_{(\lambda+\rho, \varepsilon_{i})}+(-1)^{\p{v}_{i}}\gamma_{(\lambda+\rho-\varepsilon_{i}, \varepsilon_{i})}.$  From this we conclude that \eqref{E:Eahomspace} holds if and only if $\mu=\lambda-\varepsilon_{i}$ for some $1 \leq i \leq m+n$ which is $r$-normal for $\lambda.$

Recall from \eqref{E:dualfiltration} the filtration of $W(\lambda)\otimes V^{*}$ by co-Verma supermodules.  Intersecting with $L(\lambda)\otimes V^{*}$ and projecting onto $E_{r}L(\lambda),$ we obtain a filtration 
\[
0=Q_{0} \subseteq \dots \subseteq Q_{m+n} = E_{r}L(\lambda)
\] 
of $E_{r}L(\lambda)$ where each $Q_{k}/Q_{k-1}$ is a (possibly zero) submodule of a co-Verma module.  Refining this filtration by requiring strict inclusions, we obtain 
\[
0=Q_{0} \subseteq \dots \subseteq Q_{l} = E_{r}L(\lambda)
\]
with $Q_{k}/Q_{k-1}$ a nonzero submodule of $W(\lambda - \varepsilon_{i_{k}})$ and $i_{1} < \dots < i_{l}.$

Now for any $\mu \in X(T)$ we have,
\begin{equation*}
\dim_{k} \Hom_{\Dist{G}}(M(\mu), E_{r}L(\lambda))\leq \sum_{k=1}^{l} \dim_{k} \Hom_{\Dist{G}}(M(\mu),Q_{k}/Q_{k-1}).
\end{equation*}
 
If $1 \leq t \leq m+n $ is $r$-normal for $\lambda,$ then the left hand side is nonzero for $\mu = \lambda-\varepsilon_{t}.$  However, by Lemma~\ref{L:homspaces}, the right hand side is one if $\mu = \lambda-\varepsilon_{i_{k}}$ for some $k=1, \dots ,l$ and zero otherwise.  Consequently, $t=i_{k}$ for some $1 \leq k \leq l.$  On the other hand, assume $Q_{k}/Q_{k-1}$ is nonzero.  Applying the functor $\Hom_{\Dist{G}}(M(\lambda - \varepsilon_{i_{k}}), -)$ to the short exact sequence 
\[
0\to Q_{k-1} \to Q_{k} \to Q_{k}/Q_{k-1} \to 0
\] we obtain the long exact sequence
\begin{multline*}
0 \to \Hom_{\Dist{G}}(M(\lambda-\varepsilon_{i_{k}}), Q_{k-1}) \to \Hom_{\Dist{G}}(M(\lambda -\varepsilon_{i_{k}}), Q_{k}) \to\\
 \to\Hom_{\Dist{G}}(M(\lambda - \varepsilon_{i_{k}}), Q_{k}/Q_{k-1})
\to \Ext_{\Dist{G}}^{1}(M(\lambda - \varepsilon_{i_{k}}), Q_{k-1}) \to \dots
\end{multline*}  By weights we have $\Hom_{\Dist{G}}(M(\lambda-\varepsilon_{i_{k}}), Q_{k-1})=0$ and by Lemma~\ref{L:exttriviality} we have $ \Ext_{\Dist{G}}^{1}(M(\lambda - \varepsilon_{i_{k}}), Q_{k-1})=0.$  Consequently, by Lemma~\ref{L:homspaces}, 
\begin{align*}
0 \neq \Hom_{\Dist{G}}(M(\lambda - \varepsilon_{i_{k}}), Q_{k}/Q_{k-1}) &\cong \Hom_{\Dist{G}}(M(\lambda -\varepsilon_{i_{k}}), Q_{k})\\
 &\subseteq \Hom_{\Dist{G}}(M(\lambda- \varepsilon_{i_{k}}), E_{r}L(\lambda)),
\end{align*}
hence $i_{k}$ is $r$-normal for $\lambda.$
\end{proof}

\begin{theorem}  Let $\lambda \in X(T)$ and $r \in \Z / p\Z.$
\begin{enumerate}
\item [(i)]  $E_{r}L(\lambda) \neq 0$ if and only if $\varepsilon_{r}^{*}(\lambda) \neq 0,$ in which case it is a self-dual indecomposable module with irreducible socle and cosocle both isomorphic to $L(\mu)$ where $\mu=\tilde{e}^{*}_{r}(\lambda).$  Moreover, $E_{r}L(\lambda)$ is irreducible if and only if $\varepsilon_{r}^{*}(\lambda)=1.$ 
\item [(ii)]   $F_{r}L(\lambda) \neq 0$ if and only if $\varphi_{r}^{*}(\lambda) \neq 0,$ in which case it is a self-dual indecomposable module with irreducible socle and cosocle both isomorphic to $L(\mu)$ where $\mu=\tilde{f}^{*}_{r}(\lambda).$  Moreover, $F_{r}L(\lambda)$ is irreducible if and only if $\varphi_{r}^{*}(\lambda)=1.$ 
\end{enumerate}
\end{theorem}

\begin{proof}  First we prove $(i).$  Using Corollary~\ref{C:goodcriteria} and arguing as in the proof of the previous lemma, we have
\begin{equation*}
\Hom_{\Dist{G}}(L(\mu), E_{r}L(\lambda)) \neq 0
\end{equation*}
holds if and only if $\mu=\lambda-\varepsilon_{i}$ for some $1 \leq i \leq m+n$ which is $r$-good for $\lambda.$  Also we note that if $i$ is $r$-good for $\lambda,$ then it is necessarily unique by the definition of $r$-good.

From this we see that if $E_{r}L(\lambda) \neq 0,$ then the socle of $E_{r}(\lambda)$ is precisely $L(\lambda -\varepsilon_{i}) \cong L(\tilde{e}_{r}^{*}(\lambda))$ where $i$ is $r$-good for $\lambda.$  Therefore, $E_{r}L(\lambda) \neq 0$ if and only if $\varepsilon_{r}^{*}(\lambda) \neq 0$ and, if it is nonzero, it has irreducible socle and is indecomposable.  The self-duality statement and, hence, the description of the head follows from the fact that $L(\lambda)$ is self-dual with respect to contravariant duality and $E_{r}$ commutes with this duality.  

Finally, we prove the statement about the simplicity of $E_{r}L(\lambda).$  If $E_{r}L(\lambda)$ is simple, it is clear that $i$ is $r$-normal if and only if $i$ is $r$-good.  Consequently, $\varepsilon_{r}^{*}(\lambda)=1.$  One the other hand, assume $\varepsilon_{r}^{*}(\lambda)=1$ but $E_{r}L(\lambda)$ is not irreducible.  Fix $i$ to be the position of the unique $-$ in the reduced $r$-signature for $\lambda.$  From Lemma~\ref{L:homsintoErLlambda} we deduce that $E_{r}L(\lambda)$ is a submodule of $W(\lambda - \varepsilon_{i}).$  This implies by the reducibility assumption that $\Hom_{\Dist{G}}(E_{r}L(\lambda),L(\lambda - \varepsilon_{i}))=0.$  However, taking contravariant duals and using that $i$ is $r$-good for $\lambda,$ we see that 
\[
\Hom_{\Dist{G}}(E_{r}L(\lambda),L(\lambda-\varepsilon_{i}))\cong \Hom_{\Dist{G}}(L(\lambda-\varepsilon_{i}),E_{r}L(\lambda)) \neq 0.
\] This gives the desired contradiction.  Therefore $E_{r}L(\lambda)$ is irreducible.

One deduces $(ii)$ from $(i)$ using the map $\tilde{S}$ defined in \eqref{E:Stildedef} or can be proven directly with an argument similar to the one used above.
\end{proof}

Since
\[
\begin{array}{ccc}
L(\lambda)\otimes V^{*} = \bigoplus_{r \in \Z/pZ} E_{r}L(\lambda) & \text{ and }  &  L(\lambda)\otimes V = \bigoplus_{r \in \Z/pZ} F_{r}L(\lambda)
\end{array}
\]
are decompositions into indecomposable summands we obtain the following corollaries.

\begin{corollary}\label{C:socles}  Let $\lambda \in X(T).$  Then,
\begin{enumerate}
\item[(i)] The socle of $L(\lambda)\otimes V^{*}$ is precisely
\[
\bigoplus_{\substack{1 \leq i \leq m+n\\ i \text{ is good for $\lambda$}}} L(\lambda -\varepsilon_{i}).
\]
\item[(ii)] The socle of $L(\lambda)\otimes V$ is precisely
\[
\bigoplus_{\substack{1 \leq i \leq m+n\\ i \text{ is cogood for $\lambda$}}} L(\lambda +\varepsilon_{i}).
\]
\end{enumerate}  In particular, in both cases the socle is multiplicity free and contains no more than $p$ irreducible summands.
\end{corollary}

\begin{corollary}\label{C:semisimplicity}  Let $\lambda \in X(T).$  Then,
\begin{enumerate}
\item [(i)] $L(\lambda) \otimes V^{*}$ is semisimple if and only if every $1 \leq i \leq m+n$ which is normal is good.  That is, if and only if $\varepsilon_{r}^{*}(\lambda) \leq 1$ for every $r \in \Z/p\Z$.
\item [(ii)] $L(\lambda) \otimes V$ is semisimple if and only if every $1 \leq i \leq m+n$ which is conormal is cogood.  That is, if and only if $\varphi_{r}^{*}(\lambda) \leq 1$ for every $r \in \Z/p\Z$.
\end{enumerate}
\end{corollary}

\section{Odd Reflections}\label{S:oddreflections}  The results of this article depend on our choice of a homogeneous basis for $V.$  More precisely, the crystal structure on $X(T)$ depends on the sequence of parities, $\p{v}_{1}, \dots , \p{v}_{m+n},$ of the vectors we have chosen.  In this section, we discuss how to translate from one choice to another using Serganova's odd reflections.  Note that it suffices to determine how to translate from the fixed homogeneous basis $v_{1}, \dotsc , v_{m+n}$ to the homogeneous basis $v_{1}, \dotsc , v_{i-1}, v_{i+1}, v_{i}, v_{i+2}, \dotsc , v_{m+n}$ when $\p{v}_{i}+\p{v}_{i+1}=\1$ for some $1 \leq i \leq m+n-1.$  Let $1 \leq i \leq m+n-1$ and let $\sigma_{i}:V \to V$ denote the linear map defined by sending the first basis to the second.  That is, 
\begin{equation*}
\sigma_{i}(v_{j})= v_{(i \: i+1)j},
\end{equation*} for all $1 \leq j \leq m+n$ and where $(i \: i+1) \in S_{m+n}.$  We then have an automorphism on $GL(V)$ induced by $\sigma_{i}$ which in turn induces an automorphism $s_{i}:\Dist{G} \to \Dist{G}.$  Explicitly, 
\begin{equation}\label{E:sidef}
s_{i}:e_{k,l} \mapsto e_{(i \: i+1)k,(i \: i+1)l}
\end{equation}
for all $1 \leq k,l \leq m+n.$
We can twist a $\Dist{G}$-supermodule $M$ by $s_{i}$ in the usual way, where $M^{s_{i}}=M$ as a superspace and $a.m=s_{i}(a)m$ for all $a \in \Dist{G}$ and all $m \in M.$  In particular, $L(\lambda)^{s_{i}}$ is again an irreducible $\Dist{G}$-supermodule of some highest weight.   The following lemma allows us to determine the highest weight of $L(\lambda)^{s_{i}}.$

Recall that $S_{m+n}$ acts on $X(T)$ via $x \cdot \varepsilon_{j} = \varepsilon_{xj}$ for $x \in S_{m+n}$ and  $1 \leq j \leq m+n.$  For $1 \leq i \leq m+n-1,$ let $s_{i}: X(T) \to X(T)$ be the involution given by 
\begin{align}\label{E:crystalisodef}
s_{i}(\lambda) &= \begin{cases} (i \: i+1) \cdot \lambda, & \text{ if } (\lambda, \varepsilon_{i}-\varepsilon_{i+1}) \equiv 0 \modp;\\
                               (i \: i+1) \cdot \lambda - \varepsilon_{i+1} +\varepsilon_{i}, & \text{ if } (\lambda, \varepsilon_{i}-\varepsilon_{i+1}) \not\equiv 0 \modp;
\end{cases}
\end{align} where $(-,-)$ is the bilinear form on $X(T)$ defined in \eqref{bf}.

\begin{lemma}\label{L:st}
Let $\lambda \in X(T)$ and let $1 \leq i \leq m+n-1$ satisfying $\p{v}_{i}+\p{v}_{i+1}=\1.$ 
Then, 
$$
L(\lambda)^{s_{i}} \cong L(s_{i}(\lambda))
$$
\end{lemma}

\begin{proof}

We observe that by the action of $s_{i}$ on $\Dist{T},$ if $v \in L(\lambda)_{\mu},$ then $v \in (L(\lambda)^{s_{i}})_{(i \:i+1) \cdot \mu}.$ 
Now let $v_{\lambda}$ be a $\Dist{G}$ primitive vector in $L(\lambda)$ of weight $\lambda$,
c.f. Lemma~\ref{hw}. 
We observe  that the vector $e_{i,i+1}.v_{\lambda}=e_{i+1,i} v_{\lambda}  \in L(\lambda)^{s_{i}}$ is a $\Dist{G}$ primitive vector.   For if $1 \leq r<s \leq m+n,$ then, unless $r=i$ and $s=i+1,$ it is straightforward to verify that $e_{r,s}.e_{i+1,i}v_{\lambda}=0$ by the action of $s_{i}$ and by \eqref{meq}.   The final case to check is $e_{i, i+1}.e_{i+1,i}v_{\lambda}.$  However, $ e_{i, i+1}.e_{i+1,i}v_{\lambda}=e_{i+1,i}^{2}v_{\lambda}$ and since $\p{e}_{i+1,i}=\p{v}_{i}+\p{v}_{i+1}=\1,$ $e_{i+1,i}^{2}=0$ and so the result follows in this case as well.

Now suppose that $e_{i+1,i} v_{\lambda} \neq 0$.
We get from the previous paragraph
that $e_{i+1,i} v \in  L(\lambda)^{s_{i}}$ is $\Dist{G}$ primitive of weight $(i \: i+1) \cdot (\lambda - \varepsilon_{i}+\varepsilon_{i+1})= (i \: i+1) \cdot \lambda -\varepsilon_{i+1}+\varepsilon_{i}$.
Hence, $L(\lambda)^{s_{i}} \cong L( (i \: i+1) \cdot \lambda -\varepsilon_{i+1}+\varepsilon_{i})$.
On the other hand, if $e_{i+1,i} v_{\lambda} = 0$, then $v_{\lambda} \in L(\lambda)^{s_{i}}$ is itself 
already a $\Dist{G}$ primitive vector of weight $(i \: i+1) \cdot \lambda$ so
$L(\lambda)^{s_{i}} \cong L((i \: i+1) \cdot \lambda)$.

Thus, to complete the proof of the lemma, it suffices to show that
$e_{i+1,i} v_{\lambda} \neq 0$ if and only if 
$(\lambda, \varepsilon_{i}-\varepsilon_{i+1}) \not\equiv 0 \modp$. 
But $e_{i+1,i} v_{\lambda} \neq 0$ if and only if 
there is some element $x\in\Dist{B}$ such that 
$x e_{i+1,i} v_{\lambda}$ is a non-zero multiple of $v_{\lambda}$.
By weights, the only $x$ that needs to be considered
is $e_{i,i+1}$. Finally, $e_{i,i+1} e_{i+1,i} v_{\lambda} = (-1)^{\p{v}_{i}}(\lambda,\varepsilon_{i}-\varepsilon_{i+1}) v_{\lambda}$.
\end{proof}

This result is closely related to the fact that $G$ has non-conjugate Borel subgroups and, hence, different labelings of the irreducible $\Dist{G}$-supermodules by highest weight.  The problem of translating between labelings was first solved by Serganova \cite{sergthesis} and is also discussed in \cite{BK0}.

\begin{corollary}\label{C:combinatorialRmatrices} Let $ 1 \leq i \leq m+n$ with $\p{v}_{i}+\p{v}_{i+1}=\1.$ Let $X(T)_{1}$ denote $X(T)$ with the crystal structure lifted (as in subsection~\ref{SS:crystalofXofT}) from the crystal 
\[
\left(\mathcal{B}_{\p{v}_{1}}\otimes \dotsb \otimes \mathcal{B}_{\p{v}_{i}} \otimes \mathcal{B}_{\p{v}_{i+1}}\otimes \dotsb \otimes \mathcal{B}_{\p{v}_{m+n}}, \tilde{e}_{r}^{*}, \tilde{f}_{r}^{*}, \varepsilon_{r}^{*}, \varphi_{r}^{*}, \wt  \right),
\] and let $X(T)_{2}$ denote $X(T)$ with the crystal structure lifted from the crystal 
\[
\left(\mathcal{B}_{\p{v}_{1}}\otimes \dotsb \otimes \mathcal{B}_{\p{v}_{i+1}} \otimes \mathcal{B}_{\p{v}_{i}}\otimes \dotsb \otimes \mathcal{B}_{\p{v}_{m+n}}, \tilde{e}_{r}^{*}, \tilde{f}_{r}^{*}, \varepsilon_{r}^{*}, \varphi_{r}^{*}, \wt  \right).
\] The map 
\[
s_{i}: X(T)_{1} \to X(T)_{2}
\] defined in \eqref{E:sidef} gives an isomorphism of crystals.
\end{corollary}

\begin{proof} We first prove that $\wt (s_{i}(\lambda))=\wt (\lambda).$ However, in both cases (when  $(\lambda, \varepsilon_{i}-\varepsilon_{i+1}) \equiv 0 \modp$ and when $(\lambda, \varepsilon_{i}-\varepsilon_{i+1}) \not\equiv 0 \modp$) this follows from straightforward arguments using the definition of the function $\wt$ and a careful consideration of the parities involved.

We next prove that $s_{i}(\tilde{e}_{r}^{*}(\lambda))=\tilde{e}_{r}^{*}(s_{i}(\lambda))$ for all $\lambda \in X(T)$ and all $ r \in \Z/p\Z.$ We first prove that $\tilde{e}_{r}^{*}(\lambda) \neq 0$ if and only if $\tilde{e}_{r}^{*}(s_{i}(\lambda)) \neq 0.$ If $\tilde{e}_{r}^{*}(\lambda) \neq 0,$ then by Theorem~\ref{T:transfunctorsonirreps} $\wt (\tilde{e}_{r}^{*}(\lambda))=\wt (\lambda)+\gamma_{r}-\gamma_{r+1}$ and, twisting by $s_{i}$ and applying Lemma~\ref{L:st}, 
\begin{align*}
0  & \neq \Hom_{\Dist{G}}(L(\tilde{e}_{r}^{*}(\lambda)), L(\lambda)\otimes V^{*})  \\
  & \cong \Hom_{\Dist{G}}(L(\tilde{e}_{r}^{*}(\lambda))^{s_{i}}, L(\lambda)^{s_{i}}\otimes (V^{*})^{s_{i}})\notag \\
  & \cong \Hom_{\Dist{G}}(L(s_{i}(\tilde{e}_{r}^{*}(\lambda))), L(s_{i}(\lambda)) \otimes V^{*}). \notag
\end{align*}
By the previous paragraph we have $\wt (s_{i}(\tilde{e}_{r}^{*}(\lambda)))=\wt (s_{i}(\lambda))+\gamma_{r}-\gamma_{r+1}.$  By Theorem~\ref{T:transfunctorsonirreps} again, we have
\[
 \Hom_{\Dist{G}}(L(\tilde{e}_{r}^{*}(s_{i}(\lambda))), L(s_{i}(\lambda)) \otimes V^{*}) \neq 0.
\]  Thus, $\tilde{e}_{r}^{*}(s_{i}(\lambda)) \neq 0.$  The converse is proven by in an identical manner.
We next observe that if $\tilde{e}_{r}^{*}(s_{i}(\lambda)) \neq 0,$ then it is completely determined by the fact that it satisfies the following two properties:
\begin{enumerate}
\item $\Hom_{\Dist{G}}(L(\tilde{e}_{r}^{*}(s_{i}(\lambda))), L(s_{i}(\lambda))\otimes V^{*}) \neq 0,$
\item $\wt (\tilde{e}_{r}^{*}(s_{i}(\lambda)))=\wt (s_{i}(\lambda))+\gamma_{r}-\gamma_{r+1}.$
\end{enumerate}  Consequently, by our earlier observations we have $s_{i}(\tilde{e}_{r}^{*}(\lambda))=\tilde{e}_{r}^{*}(s_{i}(\lambda)).$

A similar argument proves $s_{i}(\tilde{f}_{r}^{*}(\lambda))=\tilde{f}_{r}^{*}(s_{i}(\lambda))$ for all $\lambda \in X(T)$ and all $ r \in \Z/p\Z.$  This proves the desired result.
\end{proof}
One can also give a purely combinatorial proof of this result.

\end{document}